\documentclass[12pt]{article}
\usepackage{latexsym}
\usepackage{amssymb}
\usepackage{amsthm,amsmath}

\makeatletter
\@addtoreset{equation}{section}

\ifx\reset@font\undefined
      \let\reset@font=\relax
\fi

\title{Quantitative estimates of the convergence of the empirical
  covariance matrix in Log-concave Ensembles
}

\author{
Rados{\l}aw   Adamczak${}^{1}$ \and  Alexander E. Litvak \and Alain Pajor
       \and
Nicole Tomczak-Jaegermann${}^{2}$ }

\newcommand\address{\noindent\leavevmode%
\noindent
Rados{\l}aw  Adamczak, \\
Institute of Mathematics, \\
University of Warsaw, \\
Banacha 2, 02-097 Warszawa, Poland\\
 \texttt{\small
 radamcz@mimuw.edu.pl}

\medskip
\noindent
Alexander E. Litvak, \\
Dept.~of Math.~and Stat.~Sciences,\\
University of Alberta, \\
Edmonton, Alberta, Canada, T6G 2G1.\\
\texttt{\small%
  alexandr@math.ualberta.ca}

\medskip
\noindent
Alain  Pajor, \\
Universit\'{e} Paris-Est\\
\'{E}quipe d'Analyse et Math\'{e}matiques Appliqu\'ees, \\
5, boulevard Descartes,
Champs sur Marne,\\
77454 Marne-la-Vall\'{e}e,  Cedex 2, France\\
\texttt{\small
e-mail: Alain.Pajor@univ-mlv.fr }

\medskip
\noindent
Nicole  Tomczak-Jaegermann, \\
Dept.~of Math.~and Stat.~Sciences,\\
University of Alberta, \\
Edmonton, Alberta, Canada, T6G 2G1.\\
\texttt{\small
e-mail:    nicole.tomczak@ualberta.ca}
}

\date{}

\newtheorem{fact}{Fact}[section]
\newtheorem{thm}[fact]{Theorem}

\newtheorem{prop}[fact]{Proposition}
\newtheorem{lemma}[fact]{Lemma}
\newtheorem{cor}[fact]{Corollary}
\newtheorem{remark}[fact]{Remark}

\newtheorem{definition}[fact]{Definition}
\newbox\nrmbox
\setbox\nrmbox=\hbox{$\Vert$}
\def\nrmrule{\vrule height\ht\nrmbox depth1.2\dp\nrmbox}

\setbox\nrmbox=
\hbox{\kern0.15em\nrmrule\kern0.15em\nrmrule
   \kern0.15em\nrmrule\kern0.15em}
\newcommand{\Snorm}[1]
      {\copy\nrmbox#1\copy\nrmbox\kern-0.03em
            \lower.4ex\hbox{}}

\renewcommand{\proof}{\noindent{\bf Proof{\ \ }}}

\renewcommand{\qed}{\bigskip\hfill\(\Box\)}
\newcommand*{\ind}[1]{\mathbf{1}_{\{#1\}}}

\newcommand{\R}{\mathbb R}

\newcommand{\E}{\mathbb E}
\newcommand{\PP}{\mathbb P}
\newcommand{\Rn}{\R^n}
\newcommand{\pp}{\mbox{{\it I\kern -0.25emP}}}

\newcommand{\ep}{\varepsilon}
\newcommand{\eps}{\varepsilon}

\def\la{\left\langle}
\def\ra{\r\rangle}

\newcommand{\Id}{\mathop{\rm Id}}

\newcommand{\supp}{\mathop{\rm supp}}

\def\r{\right}

\begin{document}

\maketitle
\footnotetext[1]{Work on this paper began when this author held a
  postdoctoral position at the Department of Mathematical and
  Statistical Sciences, University of Alberta in Edmonton,
  Alberta. The position was partially sponsored by the Pacific Institute for the
  Mathematical Sciences.}
\footnotetext[2]{This author holds the Canada Research Chair in
  Geometric Analysis.}

\begin{abstract}
 Let $K$ be an isotropic convex body in $\R^n$.  Given $\eps>0$, how
many independent points $X_i$ uniformly distributed on $K$ are needed
for the empirical covariance matrix to approximate the identity up to
$\eps$ with overwhelming probability? Our paper answers this question
from \cite{KLS}. More precisely, let $X\in\R^n$ be a centered random
vector with a log-concave distribution and with the identity as
covariance matrix.  An example of such a vector $X$ is a random point
in an isotropic convex body. We show that for any $\eps>0$, there
exists $C(\eps)>0$, such that if $N\sim C(\eps)\, n$ and $(X_i)_{i\le
N}$ are i.i.d. copies of $X$, then $ \Big\|\frac{1}{N}\sum_{i=1}^N
X_i\otimes X_i - \Id\Big\| \le \varepsilon, $ with probability larger
than $1-\exp(-c\sqrt n)$.\\

\noindent AMS Classification:  primary 52A20, 46B09,  52A21  secondary
15A52, 60E15

\noindent Keywords: convex bodies, log-concave measures, isotropic
measures, random matrices, norm of random matrices, uniform laws of
large numbers, approximation of covariance matrices

\end{abstract}

\section{Introduction}
\label{intro}

Let $X\in\R^n$ be a centered random vector with covariance matrix
$\Sigma$ and consider $N$ independent random vectors $(X_i)_{i\le
N}$ distributed as $X$.  By the law of large numbers, the empirical
covariance matrix $\frac{1}{N}\sum_{i=1}^N X_i\otimes X_i$ converges
to  $\E\, X\otimes X= \Sigma$ as $N \to \infty$. Our aim is to give
quantitative
estimate of the rate of this convergence,  that is, to estimate the
size $N$ of the sample for which
\begin{equation}\label{discrepancy}
\Big\|\frac{1}{N}\sum_{i=1}^N
X_i\otimes X_i - \Sigma\Big\| \le \varepsilon\|\Sigma\|
\end{equation}
holds with high probability.

This question was investigated in \cite{KLS} motivated by a problem
of complexity in computing volume in high dimension. In particular
the authors proved that
$$
  \E \,\Big\|\frac{1}{N}\sum_{i=1}^N  X_i\otimes X_i -
  \Sigma\Big\|\le C\,{\frac{n^2}{ N} }\|\Sigma\|,
$$
where $C=\max_{i\le N} \E|X_i|^4 / (\E|X_i|^2)^{2}$.
Chebyshev's inequality yields then a first estimate:
for any $\eps>0$, $\delta\in (0,1)$,
\begin{equation}
\label{question1}
\PP \,\Big(\Big\|\frac{1}{N}\sum_{i=1}^N
X_i\otimes X_i - \Sigma\Big\|\le \eps\|\Sigma\|\Big)\ge 1-\delta
\end{equation}
whenever $N\ge {\frac{C}{ \eps\delta}} n^2$.

When random vectors are standard Gaussian, the covariance matrix is
the identity and it is known (see the survey \cite{DZ}) that
(\ref{discrepancy}) holds with high probability whenever $N\ge
4n/\eps^2$.  This raises the question about the order of the best
$N$. In particular  can it be proportional to $n$, under reasonable
assumptions? More precisely,  the question in \cite{KLS} was phrased
in the following setting.

Let $K\subset\R^n$ be a convex body and let $X\in K$ be a random
point uniformly distributed on $K$.  Suppose that $X$ is centered at
0 and that the covariance matrix of $X$ is the identity of  $\R^n$.
In such a case we shall say that $X$ (or $K$) is isotropic.  Note
that any convex body with non empty interior has an affine isotropic
image. In this setting and under these assumptions, the question may
be stated as follows:

\smallskip
\noindent{\bf Question:} (\cite{KLS}) {\em Let $K$ be an isotropic convex
  body in $\R^n$.  Given $\eps>0$, how many independent points $X_i$
  uniformly distributed on $K$ are needed for the empirical covariance
  matrix to approximate the identity up to $\eps$ with overwhelming
  probability?}

Our main aim in this paper is to answer this question. As it is well
known to specialists, a good framework  for this kind of geometric
probabilistic questions is given by log-concave distribution (see below
for the  definition).  This is a stable and well structured class of
measures in $\R^n$ that contains uniform measure on convex bodies.
Thus our goal is to estimate
\begin{equation}\label{probability_discrepancy}
\PP \,\Big(\Big\|\frac{1}{N}\sum_{i=1}^N
X_i\otimes X_i - \Sigma\Big\|\le \eps\|\Sigma\|\Big)
\end{equation}
where $\Sigma$ is the covariance matrix of a centered random vector
$X\in\R^n$ with a log-concave distribution and $(X_i)$ are $N$
independent random vectors distributed as $X$.

Since for a symmetric matrix $M$, one has $\|M\| = \sup_{y\in S^{n-1}}
\langle M y, y\rangle$, (\ref{discrepancy}) is implied by
\begin{equation}
\label{equiv_discrepancy}
  \Big|\frac{1}{N}\sum_{i=1}^N(\langle X_i,y\rangle^2 -
  \E\langle X_i,y\rangle^2)\Big|\le \varepsilon \langle \Sigma
  y,y\rangle \quad \text{for
    all} \quad {y\in \R^n}.
\end{equation}
In the case when the covariance matrix is the identity, it is
equivalent to
\begin{equation}\label{isotropic case}
 1-\eps\le \frac{1}{N}\sum_{i=1}^N\langle X_i,y\rangle^2\le 1+\eps
 \quad \text{for all} \quad {y\in S^{n-1}}.
\end{equation}

Because of the linear invariance, there is no loss of generality to
consider  just this case when the covariance matrix is the identity.

In this framework, a breakthrough was achieved in \cite{B} where it
was proved that for any $\eps,\delta\in (0,1)$, there exists
$C(\eps,\delta)>0$ such that if a body $K$ is isotropic then $N=
C(\eps,\delta)n\log^3n$ i.i.d. uniformly distributed points on $K$
satisfy (\ref{question1}).  This estimate was further improved to $N=
C(\eps,\delta)n\log^2 n$ in \cite{R} and to $N= C(\eps,\delta)n\log n$
in \cite{GHT} and \cite{Pao}; the former paper treated the case when
$K$ is invariant under every reflection with respect to coordinate
subspaces and the latter proved the estimate in full generality

One should note that in all these results, the probability in
(\ref{question1}) does not go to 1 as $n$ goes to infinity, as one
expects in this type of high dimensional phenomena. This probability,
$1-\delta$, is given by a parameter $\delta$ and $C(\eps, \delta)$
depends on it. Thus letting $\delta $ tend to zero may destroy the
estimate on $N$.  To emphasize this important feature we will talk
about {\em overwhelming probability} if the probability goes to 1 as
$n$ goes to infinity.

The first result establishing (\ref{discrepancy}) with overwhelming
probability was given in \cite{MP}.  When a body $K$ is invariant
under every reflection with respect to coordinate subspaces, it is
proved in \cite{A} that for any $\eps\in (0,1)$ there exist
$C(\eps)>0$ such that (\ref{isotropic case}) holds whenever $N\ge
C(\eps)\,n$ and with probability going to 1 as $n$ goes to infinity.
Finally, the present paper shows, as a consequence of our main results
(Theorems \ref{inertia_matrix_thm} and \ref{empirical_moments_thm}),
that the same is true for an arbitrary body $K$ (in the isotropic
position).

\medskip

An important related direction concerns norms of random matrices with
independent log-concave columns (or rows).  More precisely, let
$X\in\R^n$ be a centered random vector with a log-concave distribution
such that the covariance matrix is the identity.  Consider $N$
independent random vectors $(X_i)_{i\le N}$ distributed as $X$ and
define $A=A^{(N)}$ to be the $n\times N$ matrix with $(X_i)_{i\le N}$
as columns.  For $n, N$ arbitrary (and $N$ not too large, actually, $n
= N$ being the central case) the question is to prove an estimate for
the norm $\|A\|$ as an operator $A: \ell_2^N \to \ell_2^n$, valid with
overwhelming probability. This problem can be viewed as an
``isomorphic form'' of an upper estimate in (\ref{isotropic case})
(for $n=N$, say), and the papers discussed above provided some answers
-- with ``parasitic'' logarithmic factors -- to this question as well.
The present article gives optimal estimates for $\|A\|$ (in Theorem
\ref{norm_estimate_thm} and Corollaries \ref{norm_cor} and
\ref{ell_p_cor}); for example, for the square matrix if $n=N$, we have
$\|A\|\le C \sqrt n$, with overwhelming probability.

\medskip

To observe a still one more point of view, for arbitrary $n$ and $N$,
consider again $A=A^{(N)}$.  The set of $n\times n$ matrices may be
equipped with the distribution of $AA^*$ to be a matrix
probability space and because of the analogy with Random Matrix
Theory, in particular with Wishart Ensemble, let us call it a {\em
  Log-concave Ensemble}.

In the last decades, in Asymptotic Geometric Analysis, considerable
work and progress have been achieved in understanding the properties
of random vectors with log-concave distribution, and more recently, in
understanding spectral properties of random matrices with independent
rows (or columns) with log-concave distribution.  It appears that in
high dimension they behave somewhat similarly as if the coordinate
would be independent.  This leads by analogy with Random Matrix Theory
to questions on the spectrum of $AA^*$ similar to those of the Wishart
Ensemble. One important difference is that now the entries are
dependent but strongly structured by the log-concavity hypothesis.

Denote by
$\lambda_1=\lambda_1(A^{(N)})\le\cdots\le\lambda_n=\lambda_n(A^{(N)})$
the eigenvalues of $AA^*$ (the squares of the singular values of
$A$). It was proved in \cite{PP} that when $n/N$ goes to
$\beta\in(0,1)$ as $n,N \to \infty$, then the empirical measures of
the eigenvalues have
a limit. It is the so-called Marchenko-Pastur distribution, as for
the Wishart Ensemble when all entries of the matrix $A$ are i.i.d.
It is also known (\cite{BY}) in the case when all the entries of $A$
are i.i.d. (with a finite fourth moment) and $\lim _{n\to
+\infty}{\frac{n}{ N}}=\beta\in (0,1)$ that $\lim
\lambda_1/N= (1-\sqrt{\beta})^2$ and $\lim \lambda_n/N= (1+\sqrt{\beta})^2$.
One could conjecture that such
results are also valid in the log-concave setting.  Nevertheless,
these results are asymptotic and not quantitative (given fixed
dimension).

Problem (\ref{isotropic case}) is of course equivalent to quantitative
estimates for $\lambda_1(A^{(N)})$ and $\lambda_n(A^{(N)})$, that is of
the support of the spectrum of $A$.  An answer is given by Proposition
\ref{almost_isometric_prop}
where it is shown that for $n\le N\le
\exp(\sqrt n)$,
$$
 1-C\sqrt{\frac{n}{ N}}\log{\frac{2N}{ n}}\le \frac{1}{N}\sum_{i=1}^N \langle
 X_i,y\rangle^2\le 1+
 C\sqrt{\frac{n}{ N}}\log{\frac{2N}{ n}} \quad \quad \text{for all} \ {y\in S^{n-1}}
$$
holds with probability larger than $1-\exp(-c\sqrt n)$, where $C,c>0$
are numerical constants.  Thus, putting $\beta={\frac{n}{ N}}\in(0,1)$, we
get
$$
 1-C\sqrt\beta\log{(2/\beta)}\le
 \frac{\lambda_1}{N}\le \frac{\lambda_n}{N}\le 1+
 C\sqrt{\beta}\log{(2/\beta)}
$$
with overwhelming probability.  As a consequence already mentioned earlier,
$\|A\|\le C(\sqrt N+ \sqrt n)$ with overwhelming probability, where
$C>0$ is a numerical constant (Corollary \ref{ell_p_cor}).

\bigskip

Our  general method follows an approach that can be traced back to
Bourgain \cite{B} (cf.~also \cite{GM}). It relies upon  a crucial
new ingredient of a novel chaining argument that in an essential way
depends on the distribution of coordinates of a point on the unit
sphere. What makes this approach work, by rather subtle estimates,
is a special structure of the sets used for the chaining.

To describe a very rough idea of this structure, involved in the proof
of Theorem \ref{norm_estimate_thm} below, assume for simplicity that
$m = n = 2^{s}$ and let $a_k = 2^{s-k}$ for $1\le k \le s$.  For each
$k$, first consider the subset of the Euclidean unit ball in $\R^N$ of
all vectors that have the support of cardinality less than or equal to
$a_k$ and with the $\ell_\infty$ norm of the coordinates bounded by
$\alpha_k$, and then define ${\cal M}^{(k)}$ to be a preassigned
$\eps_k$ net (in the Euclidean norm) of this set, where $ 0 <
\alpha_k,\, \eps_k <1$ are judiciously fixed in advance. Using sets
${\cal M}^{(k)}$ in successive steps of chaining we arrive to the set
$\cal M$ that consists of sums $v= \sum_k v_k$ where $v_k$'s are
mutually disjointly supported vectors from ${\cal M}^{(k)}$ (assuming
that the Euclidean norm of $v$ is less than 2).  As can be expected
the actual definition of $\cal M$ contains a number of delicate points
which were omitted here and can be found at the beginning of the proof
of Theorem \ref{norm_estimate_thm}.  However it is given in just one
step without discussing each individual step of the chaining.

\bigskip

The paper is organized as follows. In the next
 Section \ref{notation_section}
 we present some
definitions and
preliminary tools.  In
Section \ref{norm_section}
we study the norm of a restriction of
the matrix $A = A^{(N)}$ defined by
$$
     A_m = \sup _{F\subset \{1, ..., N\} \atop |F|\leq m } \| A _{|
       \R^F} \| =
     \sup _{z\in S^{N-1} \atop |\supp z| \leq m} |Az|.
$$
We show in Theorem \ref{norm_estimate_thm}
that with overwhelming
probability,
$$
A_m \leq C  \left(\sqrt{n} + \sqrt{m}   \log \frac{2N}{m} \r).
$$
In Section   \ref{cov_matrix_subsection}
we prove the result
announced in the abstract, answering a question from \cite{KLS}.  This
theorem appears as a particular case of a more general study of
$$
\sup_{y\in S^{n-1}}\Big|\frac{1}{N}\sum_{i=1}^N(\langle X_i,y\rangle^p
- \E\langle X_i,y\rangle^p)\Big|
$$
defined for any $p\ge 1$.  Such processes have been studied in
\cite{GM}, \cite{GR} and \cite{Men}.

Section \ref{additional_obs_subsection}
describes several observations
for norms of random matrices from $\ell_2$ to $\ell_p$, $p \ne 2$.  In
the final Section \ref{elementary_section}
we sketch a more elementary
proof of the main result of Section \ref{cov_matrix_subsection},
when  $p=2$.

\section{Notation and preliminaries}
\label{notation_section}

We equip $\R^n$ and $\R^N$ with the natural scalar product
$\langle\,\cdot, \,\cdot\rangle$ and the natural Euclidean norm
$|\cdot|$. We also denote by the same notation $|\cdot|$ the
cardinality of a set.
In this paper, $X$ will denote a random vector in $\R^n$
and $(X_i)$ will be independent random vectors with the same
distribution as $X$. By $\Id$ we shall denote the identity on $\R^n$
and by $\Sigma=\Sigma(X)=\E\, X\otimes X$, the covariance matrix of
$X$ (here $X\otimes X$ is the rank one operator defined by $X\otimes X
(y)= \langle X, y\rangle X$, for all $y\in\R^n$).  By $\|M\|$ we shall
denote the operator norm of a matrix $M$, that is
$\|M\|=\sup_{|y|=1}|My|$.

\begin{definition}
\label{isotropic_vector}
A random vector $X\in \R^{n}$ is called isotropic if
\begin{equation}
\E\langle X,y\rangle=0,\quad \E\,|\langle X, y \rangle|^{2}=|y|^{2}
\quad \mbox{\rm for all }
 y\in \R^{n},  \label{isotropic}
\end{equation}
in other words, if $X$ is centered and its covariance matrix is the identity:
$$
\E\, X\otimes X=  \Id.
$$
\end {definition}

Recall that a function $f: \R^n \to \R $ is
called log-concave if for any $\theta \in \lbrack 0,1]$ and any
$x_{1},x_{2}\in \mathbb{R}^{n}$,
\begin{equation*}
f\big(\theta x_{1}+(1-\theta )x_{2}\big)\geq f(x_{1})^{\theta
}f(x_{2})^{1-\theta }.
\end{equation*}

\begin{definition}
\label{logconcave_measure}
  A measure $\mu $ on $\mathbb{R}^{n}$ is log-concave if for any
  measurable subsets $A,B$ of \thinspace\ $\mathbb{R}^{n}$ and any
  $\theta \in \lbrack 0,1]$,
\begin{equation*}
\mu (\theta A+(1-\theta )B)\geq \mu (A)^{\theta }\mu (B)^{(1-\theta )}
\end{equation*}
whenever the set
\begin{equation*}
\theta A+(1-\theta )B=\{\theta x_{1}+(1-\theta )x_{2}\,:\,x_{1}\in
A,\;x_{2}\in B\}
\end{equation*}
is measurable.
\end{definition}

The Brunn-Minkowski inequality
provides examples of log-concave measures, that are the uniform Lebesgue
measure on compact convex subsets of $\mathbb{R}^{n}$ as well as their
marginals
(cf. e.g., \cite{Schneider}).
More generally,  Borell's theorem
\cite{Bor} characterizes the log-concave measures that are not
supported by any hyperplane as the absolutely continuous measures
(with respect to the Lebesgue measure) with a log-concave density. Note
that the distribution of an isotropic vector is not supported by any
hyperplane. Moreover, it is known \cite{Bor2} that if a measure is
log-concave then linear functionals exhibit a sub-exponential decay.
To be more precise, recall that for a random variable $Y$, the
$\psi_1$ norm of $Y$ is
$$
\|Y\|_{\psi_1}=\inf\left\{C>0\,;\,
\,\E\exp\left(\frac{|Y|}{C}\right) \leq 2\right\}.
$$

A straightforward computation shows that for every integer $p \geq 1$,
\begin{equation}\label{psi_norm}
(\E|Y|^p)^{1/p} \leq cp \| Y\|_{\psi_1}
\end{equation}
where $c$ is an absolute constant.

We can now state the sub-exponential decay of linear functionals in terms of
$\psi_1$ norm \cite{Bor2}:

\begin{lemma}\label{psi}
  Let $X\in\R^n$ be a centered random vector with a log-concave
  distribution.  Then for every $y \in S^{n-1}$,
$$
 \|\la X, y \ra \|_{\psi_1} \leq  \psi \,(\E|\langle X, y \rangle |^2)^{1/2}
$$
where $\psi>0$ is universal constant.
Moreover, if $X$ has a symmetric distribution then
$\psi =2$.
 \end{lemma}

The moreover part easily follows by a direct calculation
(see \cite{MiPa}).

Putting together (\ref{psi_norm}) and Lemma \ref{psi}, we get that for
every $y\in S^{n-1}$,
\begin{equation}\label{psi1}
(\E|\la X, y \ra |^p)^{1/p} \leq Cp \,(\E|\langle X, y \rangle |^2)^{1/2}
\end{equation}
where $C$ is an absolute positive constant.

\section{Norm of a random matrix}
\label{norm_section}

In this Section $X_1,\ldots,X_N$ are independent random vectors in $\R^n$.
Mostly we work with i.i.d. random vectors, distributed according
to an isotropic, log-concave probability measure on $\R^n$.
Random $n\times N$ matrix whose columns are $X_i$'s is denoted by $A$
and its operator norm from $\ell _2^N$ to $\ell _2^n$ is denoted by
$\|A\|$. We will also use the following related notation, for $1 \le m
\le N$,
$$
     A_m = \sup _{F\subset \{1, ..., N\} \atop |F|\leq m } \| A _{| \R^F} \| =
     \sup _{z\in S^{N-1} \atop |\supp z| \leq m} |Az|.
$$
Note that $A_m$ is increasing in $m$.
Given a set $E\subset \{1, ..., N\}$ by $P_E$ we denote the orthogonal
projection from $\R^N$ onto coordinate subspace of vectors whose
support is in $E$. Such a subspace is denoted by $\R^E$.

\begin{lemma}\label{max}
Let $X_1,\ldots,X_N$ be i.i.d. random vectors, distributed according
to an isotropic, log-concave probability measure on $\R^n$.
  There exists an absolute positive constant $C_0$ such that for any
  $N\leq \exp(\sqrt{n})$ and for every $K\geq 1$ one has
$$
   \max _{i\leq N} |X_i| \leq C_0 K \sqrt{n}
$$
with probability at least $1 - \exp(-K \sqrt{n})$.
\end{lemma}

\proof
By \cite{Pao} we have for every $i\leq N$
$$
  \PP \left\{ |X_i| \geq C t \sqrt{n}  \r\} \leq \exp(-t c \sqrt{n}) ,
$$
where $C$ and $c$ are absolute positive constants. The result follows
by the union bound (and adjusting absolute constants).
\qed

\begin{lemma}\label{dec}
Let $x_1, \ldots, x_N \in \Rn$. There exists a set
$E\subset \{1, ..., N\}$, such that
$$
  \sum_{i\ne j} \la x_i, x_j \ra \leq 4 \sum_{i\in E}
  \sum_{j\in E^c} \la x_i, x_j \ra .
$$
\end{lemma}

\proof
Clearly one has
$$
  2^{N-2} \sum_{i\ne j} \la x_i, x_j \ra  =
  \sum_{E\subset \{1, ..., N\} }\sum_{i\in E}
  \sum_{j\in E^c} \la x_i, x_j \ra
  \leq 2^N \max_{E\subset \{1, ..., N\} }
  \sum_{i\in E}  \sum_{j\in E^c} \la x_i, x_j \ra
$$
from which the lemma follows.
\qed

\bigskip

Now, given  a  $E\subset \{1, ..., N\}$, $\eps, \alpha \in (0, 1]$,  by
${\cal{N}}(E, \eps, \alpha)$ we denote an $\eps$-net of
$B_2^N \cap \alpha B_{\infty}^N \cap \R^E$ in the Euclidean metric.
Standard volume estimate shows that we may assume that the cardinality of
${\cal{N}}(E, \eps, \alpha)$ does not exceed
$(3/\eps)^m$, where $m$ is the cardinality of $E$.

We will need the following two  lemmas.

\begin{lemma}
   \label{mainl}
Let  $X_1,\ldots,X_N$ be independent random vectors in $\R^n$
and let $\psi >0$
such that
$$
   \sup _{i\leq N} \ \sup _{y\in S^{n-1}} \| \la X_i , y \ra \|
    _{\psi _1}  \leq \psi .
$$
   Let $m\leq N$, $\eps, \alpha \in (0, 1]$ and $L\geq 2 m \log
   \frac{12 e N}{m \eps}$. Then
$$
\PP \left(\sup _{F\subset \{1, ..., N\} \atop |F|\leq m } \ \sup
  _{E\subset F} \ \ \sup _{z\in {\cal{N}}(F, \eps, \alpha)} \
  \sum_{i\in E} \left| \la z_i X_i, \sum_{j\in F\setminus E} z_j X_j
    \ra \r| > \psi \, \alpha L A_m \r) \leq e^{-L/2}.
$$
\end{lemma}

\medskip

\proof
Denote the underlying probability space by $\Omega$.
For
$F\subset \{1, ..., N\}$ with $|F|\leq m$, $E\subset F$,
and $z\in {\cal{N}}(F, \eps, \alpha)$,
define the subset $\Omega(F, E, z)$ of $\Omega$ by
$$
\Omega(F, E, z) = \left\{\sum_{i\in E} \left| \la z_i X_i, \sum_{j\in
      F\setminus E} z_j X_j \ra \r| > \psi \alpha L A_m \right\}.
$$
Fix $F$, $ E $ and $z$ as above and
 set $y= \sum_{j\in F\setminus
  E} z_j X_j$.  Clearly, $y$ is independent of  vectors
$X_i$'s, $i\in E$, and $|y| \leq A_m$.
Note that   $|y| >0$ on  $ \Omega(F, E, z)$ (otherwise
$\la z_i X_i, y \ra =0 $ for all $i \in E$ and the sharp
inequality defining
$\Omega (F, E, z)$ would be violated).
Thus, using the fact that
$\|z\|_{\infty}\leq \alpha$, we obtain
$$
\sum_{i\in E} \left| \la z_i X_i, \sum_{j\in F\setminus E} z_j X_j \ra
  \r| \leq \alpha A_m \sum_{i\in E} \left| \la X_i, y/|y| \ra \r| ,
$$
on $\Omega (F, E, z)$.
Since $A_m  > 0$ on $\Omega (F, E, z)$,  this implies
$$
\PP\left( \Omega(F, E, z) \right)
\le
  \PP \left( \sum_{i\in E} \left| \la X_i, y/|y| \ra \r| > \psi L
      \r).
$$
On the other hand, by Chebyshev's inequality and the assumption on the
$\psi_1$-norms of linear functionals, the latter probability is less than
$$
 e^{-L} \
  \E \exp\left( \sum_{i\in E} \frac{\left| \la X_i, y/|y| \ra \r|}{\psi} \r)
  \leq  2^{|E|}\ e^{-L} \leq  2^m \ e^{-L}.
$$
Therefore by  the union bound,
\begin{eqnarray*}
  \lefteqn{
  \PP \left(\sup _{F\subset \{1, \ldots, N\} \atop |F|\leq m } \  \sup
    _{E\subset F} \ \
  \sup _{z\in {\cal{N}}(F, \eps, \alpha)} \ \sum_{i\in E} \left| \la z_i X_i,
  \sum_{j\in F\setminus E} z_j X_j \ra \r| > \psi\,\alpha L A_m \r)} \\
  &  \leq &  \sum_{k=1}^m {N \choose k} \ \ 2^m \ \left( \frac{3}{\eps}
    \r)^m \sup_{F, E, z} \PP \left(\Omega(F, E, z)\right) \\
  &  \leq &  \sum_{k=1}^m {N \choose k} \ \ 2^m \ \left( \frac{3}{\eps}
    \r)^m \ 2^m \ e^{-L} \leq
   \left(\frac{e N}{m}\r) ^m   \ \left( \frac{12}{\eps} \r)^m  \ e^{-L}\\
   &  = & \exp\left( m \log \frac{12 e N}{m \eps}  -L \r),
\end{eqnarray*}
which implies the result.
\qed

We will also need another lemma of a similar type. We provide the
proof for  sake of completeness.

\begin{lemma}
\label{ltwo}
Let  $X_1,\ldots,X_N$ be independent random vectors in $\R^n$
and let $\psi >0$
such that
$$
   \sup _{i\leq N} \ \sup _{y\in S^{n-1}} \| \la X_i , y \ra \|
    _{\psi _1}  \leq \psi .
$$
Let $1\leq k,  m\leq N$, $\eps, \alpha \in (0, 1]$, $\beta >0$, and $L>0$.
Let $B(m, \beta)$ denote the set of vectors $x\in \beta B_2^N$ with
$|\supp  x| \leq m$ and let $\cal{B}$ be a subset of $B(m, \beta)$ of
cardinality $M$. Then
\begin{eqnarray*}
\PP \left(\sup _{F\subset \{1, \ldots, N\} \atop |F|\leq k }\ \sup
  _{x\in {\cal{B}}} \ \ \sup _{z\in {\cal{N}}(F, \eps, \alpha)} \right.
 & &\!\!\!\!\!\! \left.  \sum_{i\in F} \left| \la z_i X_i,
     \sum_{j\not\in F}    x_j X_j \ra \r|
    >  \psi \alpha \beta L A_m \r)  \\
  &\leq &  M \left( \frac{6 e N}{k \eps} \r)^k    e^{-L} .
\end{eqnarray*}
\end{lemma}

\medskip

\proof
The proof is analogous to the argument
in Lemma~\ref{mainl}.
For $F\subset \{1, ..., N\}$ with $|F|\leq k$, $x\in {\cal{B}}$, and
$z\in {\cal{N}}(F, \eps, \alpha)$
consider
$$
\Omega(F, x, z) =
\left\{  \sum_{i\in F} \left| \la z_i X_i,
     \sum_{j\not\in F}    x_j X_j \ra \r|
    > \psi \alpha \beta L A_m \r\}.
$$

Fix $F$, $x$, $z$ as above and set $y= \sum_{j\not\in F} x_j X_j$.
Clearly, $y$ is independent of the vectors $X_i$'s, $i\in F$,
moreover, $|y| \leq \beta A_m$, and, similarly as
in before, $|y| >0$ on
$\Omega(F, x, z) $.  Thus, using the fact that $\|z\|_{\infty}\leq
\alpha$, we obtain
$$
  \sum_{i\in F} \left| \la z_i X_i, \sum_{j\not\in F} x_j X_j \ra \r| \leq
  \alpha \beta A_m \sum_{i\in F} \left| \la X_i, y/|y| \ra \r| ,
$$
on $\Omega(F, x, z)$.  Therefore,  again  as in Lemma~\ref{mainl},  we have
\begin{eqnarray*}
\PP\left( \Omega(F, x, z) \right)
&\le &
    \PP \left( \sum_{i\in F} \left| \la X_i, y/|y| \ra \r| > \psi L
        \r) \\
      &\leq &  e^{-L} \
  \E \exp\left( \sum_{i\in F} \frac{\left| \la X_i, y/|y| \ra
        \r|}{\psi} \r)
     \leq   2^{|F|}\ e^{-L} \leq  2^k \ e^{-L}   .
\end{eqnarray*}
By the union bound we get
\begin{eqnarray*}
  \lefteqn{
\PP \left(\sup _{F\subset \{1, ..., N\} \atop |F|\leq k }\ \sup _{x\in
    {\cal{B}}} \ \ \sup _{z\in {\cal{N}}(F, \eps, \alpha)} \
  \sum_{i\in F} \left| \la z_i X_i, \sum_{j\not\in F} x_j X_j \ra \r|
    > \psi \alpha \beta L A_m \r) }\\
   &  \leq &  M \sum_{l=1}^k {N \choose l}  \ \left( \frac{3}{\eps}
     \r)^k \ 2^k \ e^{-L} \leq
   M \left(\frac{e N}{k}\r) ^k   \ \left( \frac{6}{\eps} \r)^k  \ e^{-L},
\end{eqnarray*}
which proves the result.
\qed

\medskip

\begin{remark}
{\rm
Observe that if $X_i$'s are i.i.d. random vectors,
distributed according to an isotropic, log-concave probability measure
on $\R^n$, then, by Lemma~\ref{psi}, they satisfy the condition for
the $\psi_1$-norm of Lemmas \ref{mainl} and \ref{ltwo}.
}
\end{remark}

\begin{thm}
\label{norm_estimate_thm}
Let $n\geq 1$ and $1\leq N\leq e^{\sqrt{n}}$ be integers.
Let  $X_1,\ldots,X_N$ are  i.i.d. random vectors, distributed according
to an isotropic, log-concave probability measure on $\R^n$.
 Let $K\geq 1$.
Then there are absolute positive constants $C$ and $c$ such that
$$
  \PP \left(\exists m\leq N \, : \,  A_m \geq  C K \left(\sqrt{n} + \sqrt{m}
  \log \frac{2N}{m}  \r)  \r)  \leq \exp\left( - c K \sqrt{n} \r) .
$$
\end{thm}

\begin{remark}
\label{rem_optim}
{\rm
Let $X\in\R^n$ be a random vector with an isotropic exponential
distribution, that is with the density defined for $x=(x_i)\in\R^n$ by
$\prod_1^n{\frac{1} {\sqrt 2}} \exp(-\sqrt 2 |x_i|)$. It is clearly an
isotropic vector with a log-concave distribution. Consider now the
matrix $A^{(N)}$ build as before from a sample of $X$ of size
$N$. Since
$$
  \PP (|X|\ge t\sqrt n)\ge \int_{|s|\ge t\sqrt n} {\frac{1}
{ \sqrt 2}}
  \exp(-\sqrt 2 |s|)\, ds= \exp(-\sqrt 2 t\sqrt n)
$$
we get that for any $1\le m\le N$,
$$
   \PP( A_m\ge t\sqrt n)\ge \exp(-\sqrt 2 t\sqrt n).
$$
This shows that the probability estimate in Theorem~\ref{norm_estimate_thm}
is optimal up to numerical constants. The analysis of this example shows that
up to  numerical constants the logarithmic term in the estimate of $A_m$ in
Theorem~\ref{norm_estimate_thm} is also optimal (for the details see \cite{ALPT}).
}
\end{remark}

\medskip

Letting $m=N$ we get a clearly optimal estimate for the operator norm
$\|A\|$, valid with overwhelming probability.

\begin{cor}\label{norm_cor}
In the setting of Theorem \ref{norm_estimate_thm}  we get, for every
$K \ge 1$,
\begin{equation}
  \label{norm_bound_cor}
\|A\| \le C K \left(\sqrt{n} + \sqrt{N}  \r),
  \end{equation}
with probability at least $ 1- e^{-cK\sqrt n}$, where $C, c >0$ are
absolute constants.
\end{cor}

\begin{remark}
\label{rem_special}
{\rm
 The final remark of \cite{B} states
that by refining a bit the method of proof of Lemma 2  of that paper one may
obtain that if $X_1, \ldots, X_n$ are $n$ independent vectors in
$\R^n$ distributed according to a probability measure $\mu$ on $\R^n$
satisfying $\|\langle x, y\rangle \|_{\psi_1} < 1/\sqrt n$ for all
$y\in S^{n-1}$, then,  with probability $1-\delta$, the matrix $A$
admits the bound for the operator  norm
$$
\|A\|\le C(\delta) \left( \int \left( \max_{1\le i \le n} |X_i|\right)
\, d\mu +1 \right).
$$
By Lemmas \ref{psi} and \ref{max}, and taking into account the
normalization, this would imply a version of (\ref{norm_bound_cor})
with $N=n$  and probability $1-\delta$.
}
\end{remark}

\medskip

\begin{remark}
\label{rem_m-to-n}
{\rm
 Note that $\sqrt{n} + \sqrt{m} \log \frac{2N}{m}$ in the
formula
in Theorem~\ref{norm_estimate_thm}
 can be substituted with
$$
  \sqrt{n} + \sqrt{m} \log \frac{2N}{\max\{n, m\}}.
$$
Indeed, if $m \ge n$ there is nothing to prove, otherwise
$$
\sqrt{n} + \sqrt{m} \log \frac{2N}{m} = \sqrt{n} + \sqrt{m} \log
\frac{n}{m} + \sqrt{m} \log \frac{2N}{n} \leq 2 \sqrt{n} + \sqrt{m} \log
\frac{2N}{n}.
$$
}
\end{remark}

\medskip

Finally,  another  immediate consequence.

\begin{cor}\label{equal_weights}
  There are absolute positive constants $C$ and $c$ such that for
  every $n\geq 1$, $1\leq N \leq e^{\sqrt{n}}$, $K\geq 1$, and $X_i$'s
  as in Theorem~\ref{norm_estimate_thm} one has
$$
  \PP \left( \exists_{E\subset \{1, \ldots, N\}} \left| \sum _{i\in E} X_i \r| \geq C K \left(\sqrt{n |E|} + |E|
  \log \frac{2N}{n}  \r)  \r)  \leq \exp\left( - c K \sqrt{n} \r) .
$$
\end{cor}

\proof
Given $E$ set $m=|E|$.
Consider vector $z\in S^{N-1}$ defined by $z_i = 1/\sqrt{m}$ if $i\in E$ and
$z_i=0$ otherwise. We have
$$
   \left| \sum _{i\in E} X_i \r| = \sqrt{m} |Az| \leq   \sqrt{m} A_m.
$$
Therefore Theorem~\ref{norm_estimate_thm} and  Remark~\ref{rem_optim}
imply the result.
\qed

\medskip

\noindent
{\bf Proof of Theorem~\ref{norm_estimate_thm}.}
As $N\leq e^{\sqrt{n}}$, it is easy to see, by applying the union
bound and adjusting absolute constants, that it is sufficient to prove
that for $K$ sufficiently large and every fixed $m\leq N$, one has
$$
  \PP \left( A_m \geq C K \left(\sqrt{n} + \sqrt{m}
  \log \frac{2N}{m}  \r)  \r)  \leq \exp\left( - c K \sqrt{n} \r) .
$$

We shall define a set $\cal{M}$ of vectors with a special structure
and supports less than or equal to $m$ which serves simultaneously two
purposes: we will be able to estimate with large probability $\sup_{x
  \in \cal{M}} |A x|$, and we will use $\cal{M}$ to approximate an
arbitrary vector from $B_2^N$ of support less than or equal to $m$.
Then a standard argument will lead to the required estimate for
$A_m$.

\medskip

First observe that if for a vector $x \in S^{N-1}$ there is a
simultaneous control of the size of support
and its $\ell_\infty$-norm (more precisely, $|\supp  x| \sim s$ and
$\|x\|_\infty \le s^{-1/2}$, for some $s \ge 1$) then $|Ax|$ can be
estimated, with large probability, directly by using Lemmas~\ref{dec}
and \ref{mainl} (it is also a part of the estimates below). It is
therefore natural to expect vectors from $\cal{M}$ to be sums of
(disjointly supported) vectors admitting such a simultaneous control as
above. Formally, the definition of $\cal{M}$ splits into two cases.
If
\begin{equation}
  \label{eq:1}
 m\ \log \frac{48 e N}{m} \leq  \sqrt{n},
\end{equation}
we set
$$
   {\cal{M}} = \bigcup _{E\subset \{1,\ldots N\} \atop |E|=m}
   {\cal{N}}(E, 1/4, 1).
$$
Otherwise, let
$l$ be  the smallest integer such that
\begin{equation}
  \label{eq:2}
  \frac{m}{2^l} \log \frac{48 e 2^l  N}{m}  \leq \sqrt{n},
\end{equation}
and fix positive integers $a_0, a_1,\ldots, a_l$ such that $a_k \le
m\, 2^{-k+1}$ for $1 \le k \le l$ and $a_0 \le m\,2^{-l}$, and
$\sum_{k=0}^l a_k = m$. (We shall later set $a_k:= [m\, 2^{-k+1}]-[m\,
2^{-k}]$ for $1 \le k \le l$ and $a_0:= [m\, 2^{-l}]$.)

\smallskip

Then  set
${\cal{M}}= {\cal{M}}_0\cap 2 B_2^{N}$, where ${\cal{M}}_0$
consists of all vectors of the form $x= \sum _{k=0}^{l} x_k$, where
$x_i$'s have disjoint supports and
$$
x_0 \in \bigcup _{E\subset \{1,\ldots N\} \atop |E|\leq a_0}
{\cal{N}}(E, 1/4, 1) , \, \, x_k \in \bigcup _{E\subset \{1,\ldots N\}
  \atop |E|\leq a_k} {\cal{N}}\left(E, 2^{-k}, \sqrt{\frac{2^k}{m}}\r)
  \, \mbox{ for } \, 1 \le k \le l.
$$
Note that for every vector $x\in{\cal{M}}$ we have $|\supp  x| \leq
\sum _0^l a_k = m$ and $|x|\leq 2$.

\bigskip

We shall consider the details of the case $m\log(48eN/m) > \sqrt{n}$
(the other case, when (\ref{eq:1}) holds, can be treated similarly,
actually, it is even simpler, since the construction of $\cal{M}$ is
simpler).  Fix $x\in \cal{M}$ of the form $x= \sum _{k=0}^{l} x_k $
and let $F_k$ be the support of $x_k$ (if there are more than one such
representations, we fix one of them).  Denote the coordinates of $x$
by $x(i)$, $i\leq N$, then
\begin{eqnarray}
  |A x|^2 &=& \la \sum _{i\leq N} x(i) X_i, \sum _{i\leq N} x(i) X_i \ra =
  \sum _{i\leq N} x(i)^2 |X_i|^2 + \sum _{i\ne j} \la x(i) X_i, x(j) X_j
  \ra\nonumber\\
   &\leq & 2\max _i |X_i|^2 + D_x \leq 2 \max\{ 2\max _i |X_i|^2,
   D_x\},
\label{decompos}
\end{eqnarray}
where
$$
  D_x=\sum _{i\ne j} \la x(i) X_i, x(j) X_j \ra.
$$
Note that by Lemma~\ref{max}, $\max _i |X_i| \leq C_0 K \sqrt{n}$ with
probability
larger than $1-e^{-K \sqrt{n}}$, and we would like to get
a similar estimate for  $D_x$.

\medskip

To this aim we split $D_x$ according to the structure of $x$. Namely
we let
$$
   D'_x := \sum _{k=0}^l \sum _{i, j \in F_k \atop i\ne j} \la x(i)
   X_i, x(j) X_j \ra,
$$
and
\begin{eqnarray*}
   D''_x :& = & \sum _{k=0}^l \sum _{i \in F_k \atop j\not\in F_k} \la
   x(i) X_i, x(j) X_j \ra\\
&  = & 2 \sum _{k=1}^l  \sum _{i \in F_k } \sum _{r\in G_k} \la x(i) X_i,
  \sum _{j\in F_r}   x(j) X_j \ra ,
\end{eqnarray*}
where
$G_k = \{0, k+1, k+2, \ldots , l\}$.
Note that
$$
    D_x = D'_x + D''_x .
$$

We first estimate $D'_x$.
 By Lemma~\ref{dec} we obtain that for
every $k$ there exists a subset
$\bar F_k$ of $F_k$ such that
\begin{eqnarray*}
   D'_x  & \leq & 4 \sum _{k=0}^l \sum _{i \in \bar F_k \atop j \in
     F_k \setminus \bar F_k}
   \la x(i) X_i, x(j) X_j \ra \\
  &\leq &  4  \sup _{F\subset \{1, ..., N\} \atop |F|\leq m/2^l}\
  \sup _{E\subset F} \
   \sup _{v\in {\cal{N}}\left(F, 1/4, 1 \r)}
   \  \sum _{i \in E}  \left| \la v_i X_i, \sum _{j \in F \setminus E}
     v_j X_j \ra \r| \\
& + &
   4 \sum _{k=1}^l
   \sup_{F\subset \{1, ..., N\} \atop |F|\leq 2m/2^k}\ \sup _{E\subset F} \
   \sup_{v\in {\cal{N}}(F, 2^{-k}, \sqrt{2^k/m})}
   \  \sum _{i \in E}  \left| \la v_i X_i, \sum _{j \in F \setminus E}
     v_j X_j \ra \r|.
\end{eqnarray*}
We now apply Lemma~\ref{mainl} to each summand in the sum above with
$L=2K \sqrt{n}$, $\eps = 1/4$, $\alpha =1$ for the first summand
(note that such an $L$ satisfies the condition) and with
$L=\frac{4m}{2^k} K \log\frac{12 e N 4^k}{m}$, $\eps = 2^{-k}$,
$\alpha = \sqrt{\frac{2^k}{m}}$ for $k\geq 1$. By
the union bound we obtain
\begin{eqnarray*}
     \PP \left(\vphantom{\sqrt{\frac{2^k}{m}}}
 \sup _{x\in {\cal{M}}} D'_x \right. & >  & \left. 8\psi K
       A_m \sqrt{n}
    + 2 \psi K A_m \sum _{k=1}^l \sqrt{\frac{2^k}{m}}
   \frac{8 m}{2^k} \log\frac{12 e N 4^k}{m} \right)\\
   & \leq &  \exp\left( - K \sqrt{n} \right) +  \sum _{k=1}^l
   \exp\left( - K \frac{2m}{2^k} \log\frac{12 e N 4^k}{m} \right)\\
   & \leq &  \exp\left( - K \sqrt{n} \right) +  l \exp\left( - K
     \frac{2m}{2^l}
   \log\frac{12 e N 4^l}{m}  \right) ,
\end{eqnarray*}
where $\psi$ is the absolute constant from Lemma~\ref{psi}.

Therefore, the choice of $l$ implies the following bound, with some
absolute positive constant $C$,
\begin{eqnarray*}
     \PP \left( \vphantom{\frac{2 N}{m}}
\sup _{x\in {\cal{M}}}
  D'_x \right. & > & \left. A_m  K \left( 8 \psi \sqrt{n} + C \psi
    \sqrt{m}
  \log\frac{2 N}{m}\r)   \r)\\
 & \leq &  \exp\left( - K \sqrt{n} \r)  + l \exp\left( - K \sqrt{n} \r) \leq
   (2\sqrt{n} +1) \exp\left( - K \sqrt{n} \r).
\end{eqnarray*}
(We also used the estimate $l\leq 2\sqrt{n}$,  valid when
 $m\leq N \leq e^{\sqrt{n}}$.)

\medskip

The estimate for $D''_x$ essentially follows the same lines. In a
sense it is simpler, since we don't need to apply Lemma~\ref{dec}. For
every $1\leq k \leq l$ we consider ${\cal{M}}_k = {\cal{M}}_k' \cap
2B_2^N$, where ${\cal{M}}_k'$ consists of all vectors of the form $x=
x_0 + \sum _{s=k+1}^{l} x_s$, where $x_i$'s ($i=0, k=1, \ldots, l$)
have pairwise disjoint supports and
$$
x_0 \in \bigcup _{E\subset \{1,\ldots N\} \atop |E|\leq a_0}
{\cal{N}}(E, 1/4, 1) , \, \, x_s \in \bigcup _{E\subset \{1,\ldots N\}
  \atop |E|\leq a_s} {\cal{N}}\left(E, 2^{-s}, \sqrt{\frac{2^s}{m}}\r)
  \, \mbox{ for } \, s \ge k+1 .
$$
Then ${\cal{M}}_k \subset 2B_2^N$ and
\begin{eqnarray*}
  |{\cal{M}}_k| &\leq &  12^{a_0} \prod _{s=k+1}^l \left(3\cdot
    2^s\r)^{a_s}
  {N \choose a_s} \leq 12^{a_0} \prod _{s=k+1}^l \left(\frac{3\cdot
      2^s  e N}{a_s}\r)^{a_s}   \\
  & \leq &  \exp\left( \frac{m}{2^l} \log 12   + \sum _{s=k+1}^l
    \frac{2m}{2^s} \log \frac{3e 4^s N}{2m} \r)
  \leq \exp\left( \sum _{s=k+1}^{l+1} \frac{2m}{2^s} \log \frac{3e 4^s
      N}{2m} \r) \\
 & \leq &  \exp\left( \frac{m}{2^k} \left( \log \frac{6 e 4^k N}{m}
     \sum _{s=0}^{l-k} \frac{1}{2^s}
          +   \log 4  \sum _{s=1}^{l-k} \frac{s}{2^s} \r)  \r)
   \leq \exp\left( \frac{4 m}{2^k} \log \frac{6 e 4^k N}{m}   \r) .
\end{eqnarray*}
We also observe that
\begin{eqnarray*}
  D''_x  & = & 2 \sum _{k=1}^l  \sum _{i \in F_k } \la x(i) X_i, \sum
  _{r\in G_k}  \sum _{j\in F_r}   x(j) X_j \ra \\
  & \leq & 2 \sum _{k=1}^l   \sup _{F\subset \{1, ..., N\} \atop
    |F|\leq 2m/2^k}\
  \sup _{u\in {\cal{N}}(F, 2^{-k}, \sqrt{2^k/m})} \sup _{ v\in
    {\cal{M}} _k }
  \sum_{i\in F} \left| \la u_i X_i, \sum_{j \not\in F} v_j X_j \ra \r|  .
\end{eqnarray*}
Now we apply Lemma~\ref{ltwo} to each summand with
$$
  L = L(k) = \frac{12 m}{2^k} K \log \frac{12 e 4^k N}{m},
$$
$$
 \eps =\eps _k = 2^{-k}, \, \, \alpha =\alpha _k = \sqrt{2^k/m}, \, \,
 \beta = 2,
  \, \,   {\cal{B}} = {\cal{B}} _k = {\cal{M}} _k.
$$
Using the union bound we obtain
\begin{eqnarray*}
\lefteqn{  \PP \left( D''_x  >  48 \psi  A_m K \sum _{k=1}^l
    \sqrt{\frac{2^k}{m}} \
     \frac{m}{2^k}
   \log \frac{12 e 4^k N}{m}\r)} \\
&\leq &  \sum _{k=1}^l \exp\left(\frac{4 m}{2^k} \log \frac{12 e 4^k
    N}{m} +
  \frac{2 m}{2^k} \log \frac{3 e 4^k N}{m} - K \frac{12 m}{2^k} \log
  \frac{12 e 4^k N}{m}\r) \\
 & \leq &  \sum _{k=1}^l \exp\left( - K \frac{6 m}{2^k} \log \frac{12 e 4^k
      N}{m}\r) \leq
  l \exp\left( - K \frac{6 m}{2^l} \log \frac{12 e 4^l N}{m}\r) .
\end{eqnarray*}
As in the case for $D_x'$ it follows that
$$
    \PP \left( \sup _{x\in {\cal{M}}} D''_x > 3 C \psi A_m  K  \sqrt{m}
      \log\frac{2 N}{m}\r)
  \leq  2 \sqrt{n} \exp\left( - K \sqrt{n}  \r),
$$
where $C$ is the same absolute constant as above. Since $D_x = D'_x +
D''_x$, then
\begin{equation}
  \label{d-x-estimate}
   \PP \left( \sup _{x\in {\cal{M}}}  D_x > K A_m \left( 8 \psi \sqrt{n}
       + 4 C \psi \sqrt{m} \log\frac{2 N}{m}\r) \r)
     \leq    (4\sqrt{n} + 1)  e^{-K\sqrt n} .
\end{equation}

\bigskip

Passing now to the approximation argument,
pick an arbitrary  $z \in S^{N-1}$  with $|\supp z|\leq m$.
 Define the following subsets of $\{1, \ldots, N\}$ depending on $z$.
Denote the coordinates of $z$ by $z_i$ ($i=1, \ldots, N$).
Let  $n_1, \ldots, n_N$ be such that
$|z_{n_1}|\geq |z_{n_2}| \geq \ldots \geq |z_{n_N}|$, so that
$z_{n_i} = 0$ for $i>m$ (since $|\supp z|\le m$).
If condition (\ref{eq:1}) holds we denote the support of $z$ by
$E_0$ and consider only  this $E_0$. Otherwise  we set
$$
   E_0 = \{n_i\}_{1\leq i \leq m/2^{l}}
$$
and
$$
  E_1 = \{n_i\}_{m/2 < i\leq m}, \ E_2 = \{n_i\}_{m/4 < i\leq m/2}, \
  \ldots, \ E_l = \{n_i\}_{m/2^l < i \leq  m/2^{l-1}} ,
$$
where $l$ is the smallest integer satisfying (\ref{eq:2}) (as before).
(For small values of $n$ it can happen that $E_0$ is empty, but it does
not create any difficulty in the proof below.)
Clearly, we have
$$
  a_0  := |E_0| \leq m/2^l, \quad
 a_k := |E_k|\leq m/2^k + 1 \leq m/2^{k-1} \, \,  \mbox{ for every }
 1\leq k\leq l,
$$
and $\sum _{i=0}^l a_i=m$.
Note that the numbers  $a_k$'s  do not depend on $z$, although
the sets $E_k$'s  do.  Finally, since  $z\in S^{N-1}$,   we also observe
that  for every $k\geq 1$,
$$
    \| P_{E_k} z \| _{\infty} \leq |z_{n_s}| \le \sqrt{\frac{2^k}{m}},
$$
where $s = [m /2^{k}]$.

Note that for every
$k\geq 1$ the vector $P_{E_k} z$ can be approximated by a vector from
${\cal{N}}\left(E_k, 2^{-k}, \sqrt{\frac{2^k}{m}}\r)$ and the vector
  $P_{E_0} z$ can be approximated by a vector from ${\cal{N}}(E_0,
  1/4, 1)$.
Thus there exists  $x\in {\cal{M}}$,  with a suitable representation
$x = \sum_{k=0}^l x_k$, such that
$$
   |z-x|^2 \leq \sum _{k=0}^l |P_{E_k} z - x_k |^2 \leq
   2^{-4} + \sum _{k=1}^l 2^{-2k} < 0.4 .
$$
Moreover, $x$ is chosen to have the same support as $z$, and thus
$w = z-x$  has the support $|\supp w | \le m$.

Considering all $z \in S^{N-1}$ with $|\supp z| \leq m $
it follows that
$$
    A_m = \sup _{z\in S^{N-1} \atop |\supp z| \leq m} |A z| \leq
   \sup _{x \in {\cal{M}}} |A x| + \sqrt{0.4}
\sup _{w\in S^{N-1} \atop |\supp w| \leq m} |A w|
=  \sup _{x \in {\cal{M}}} |A x| + \sqrt{0.4}  A_m,
$$
which implies
$$
A_m \le 3   \sup _{x \in {\cal{M}}} |A x|.
$$

Recall that  by (\ref{decompos})  for every $x \in {\cal M}$ we have
$$
   |A x|^2    \leq  2 \max\{ 2\max _i |X_i|^2, D_x\},
$$
so passing to the supremum
\begin{equation}
\label{endappr}
   A_m^2 \le 9   \sup _{x \in {\cal{M}}} |A x|^2
 \leq  9 \max\{ 4\max _i |X_i|^2, 2\sup _{x\in {\cal{M}}}  D_x\}.
\end{equation}
Applying Lemma~\ref{max}  and (\ref{d-x-estimate})  we get
$$
    A_m \leq K \left( 6\, C_0 + 144 \psi \r) \sqrt{n} + 72\, C \psi K \sqrt{m}
      \log\frac{2 N}{m}
$$
with probability larger than
$$
   1 - (4\sqrt{n} +2) \exp\left( - K \sqrt{n} \r) \ge
   1-  \exp\left( - c K \sqrt{n} \r),
$$
where $c$ is an absolute positive constant.
(In fact this estimate for probability requires  that $n$ is
sufficiently  large, but, as $K\geq 1$ was arbitrary, we can
adjust the constants.)  This concludes the proof.
\qed

\medskip

\begin{remark}
  {\rm
Consider now a more general situation in which $X_1,X_2,\ldots
X_N$ -- the columns of the matrix $A$ --  are still i.i.d. centered and
log-concave, but not necessarily isotropic. Then there exists an $n
\times n$ matrix $T$, such that $(X_i)_{i=1}^N$ has the same
distribution as $(T Y_i)_{i=1}^N$, where $Y_1,\ldots,Y_N$ are
isotropic log-concave random vectors in $\R^n$. For the purpose of
computing probabilities we may assume that $X_i = TY_i$. Therefore,
with probability at least $1 - \exp(-cK\sqrt{n})$, we have for all $m
\le N$,
\begin{align*}
A_m &= \sup_{y\in S^{n-1}}\sup_{z \in S^{N-1}\atop{|\supp z| \le m}}
\Big|\sum_{i=1}^N\langle X_iz_i,y\rangle\Big|
= \sup_{y\in S^{n-1}}\sup_{z \in S^{N-1}\atop{|\supp z| \le m}}
\Big|\sum_{i=1}^N\langle Y_i z_i,T^\ast y\rangle\Big|\\
&\le \|T^\ast\|C K \left(\sqrt{n} + \sqrt{m}
  \log \frac{2N}{m}\right) = C K\kappa  \left(\sqrt{n} + \sqrt{m}
  \log \frac{2N}{m}\right),
\end{align*}
where $\kappa = \|T^\ast\| = \sqrt{\|\Sigma\|}$ (note that $\Sigma =
TT^\ast$).
}
\end{remark}

\bigskip

We conclude this section with a more technical variant of
Theorem~\ref{norm_estimate_thm}.
Note that in particular
it requires weaker conditions on $X_i$'s and  does not require
any bounds on $N$.

\begin{thm}
  \label{techn}
  Let $1\leq n$ and $1\leq N$. Let $X_1,\ldots,X_N$ be independent
  random vectors in $\R^n$ such that
$$
   \sup _{i\leq N} \ \sup _{y\in S^{n-1}} \| \la X_i , y \ra \|
    _{\psi _1}  \leq \psi .
$$
Let $A$ be a random $n\times N$ matrix whose columns are $X_i$'s, and
$A_m$, $m\leq N$, is defined as before.
Then for every $1\leq m\leq N$, every $0\leq l \leq \log m$, and every
$K\geq 1$ one has
$$
  \PP \left( A_m \geq C \psi K \left(\frac{m}{2^l} \log \frac{48 e N 2^l}{m}
  + \sqrt{m} \log \frac{2N}{m} \r) + 6 \max _{i\leq N} |X_i|  \r)
$$
$$
   \leq (1+ 2 l) \exp\left( - 2 K \frac{m}{2^l} \log \frac{12 e N 2^l}{m} \r) ,
$$
where $C$ is an absolute constant. In particular, choosing $0\leq l\leq \log m$
to be the largest integer satisfying
$$
  \frac{2 m}{2^l} \log \frac{12 e N 2^l}{m} \geq  \sqrt{m} \log \frac{2N}{m}
$$
 we obtain that for every $K\geq 1$
$$
  \PP \left( A_m \geq C \psi K \sqrt{m} \log \frac{2N}{m}
  + 6 \max _{i\leq N} |X_i|  \r)
   \leq (1+ 2 \log m ) \exp\left( -  K \sqrt{m} \log \frac{2 N}{m} \r) .
$$
\end{thm}

\medskip

\begin{remark}
  {\rm
Note that from the definitions we immediately have
$$
    A_m \geq A_1 \geq \max _{i\leq N} |X_i| .
$$
}
\end{remark}

\medskip

For completeness we outline a proof of Theorem~\ref{techn}.

\medskip

\proof (Sketch.)  We proceed as in the proof of
Theorem~\ref{norm_estimate_thm}. So first we construct $\cal{M}$. If
$l=0$ we define $\cal{M}$ exactly as after formula (\ref{eq:1}),
otherwise it will be constructed in the same way as it was constructed
after formula (\ref{eq:2}) (note that now $l$ is a fixed number). Then
we estimate $D_x= D_x' +D_x''$.  As before we use Lemmas \ref{mainl}
and \ref{ltwo}.

The only difference is that for the first summand in the formula for
$D_x'$ we use Lemma~\ref{mainl} with $L= 4 K \frac{m}{2^l}\log
\frac{48 e N 2^l}{m}$ instead of $L= 2K\sqrt{n}$. It will give us that
\begin{eqnarray*}
     \PP \left( \vphantom{\frac{2 N}{m}}
\sup _{x\in {\cal{M}}}
  D'_x \right. & > & \left.
16  A_m K \psi  \frac{m}{2^l}\log \frac{48 e N 2^l}{m}
+ C A_m K \psi  \sqrt{m}  \log\frac{2 N}{m}\r)  \\
 & \leq &  \exp\left( - 2 K \frac{m}{2^l}\log \frac{48 e N 2^l}{m} \right) +
   l \exp\left( - 2 K \frac{m}{2^l}
   \log\frac{12 e N 4^l}{m} \right)
\end{eqnarray*}
and
$$
    \PP \left( \sup _{x\in {\cal{M}}} D''_x > 3 C \psi A_m  K  \sqrt{m}
      \log\frac{2 N}{m}\r)
  \leq  l \exp\left( - K \frac{6 m}{2^l} \log \frac{12 e 4^l N}{m}\r) .
$$
Thus, with another absolute positive constant $C$ we have
\begin{eqnarray*}
     \PP \left( \vphantom{\frac{2 N}{m}}
\sup _{x\in {\cal{M}}}
  D_x \right. & > & \left.
  C A_m K \psi \left( \frac{m}{2^l}\log \frac{48 e N 2^l}{m}
  + \sqrt{m} \log\frac{2 N}{m}\r) \r)  \\
 & \leq &   \left(1+ 2 l\r) \exp\left( - K \frac{2m}{2^l}
   \log\frac{12 e N 2^l}{m} \right) .
\end{eqnarray*}

Finally we apply the same approximation procedure. By (\ref{decompos})
and approximation we get formula (\ref{endappr})
$$
  A_m^2 \leq  \max\{ 36\max _i |X_i|^2, 18 \sup _{x\in {\cal M}} D_x\},
$$
which implies the result, by adjusting constants, if necessary.  The
``in particular" part of the Theorem is trivial.  \qed

\medskip

\begin{remark}
  {\rm
It is possible to extend Theorem~\ref{techn} to a $\psi _p$-setting,
similar to the one considered in  \cite{GM}. Let $p \in [1,2]$ and let
$X$ be a random vector such that for some  $\psi _p >0$ one has
$$
   \E \exp \left( \left(|\langle X,y\rangle| /\psi _p \r) ^p \right) \le 2
$$
for every $y\in S^{n-1}$. Then, adjusting Lemmas \ref{mainl}  and \ref{ltwo},
and repeating the proof of Theorem~\ref{techn} we can get
\begin{eqnarray*}
  \lefteqn{  \PP \left( A_m \geq C \psi _p K \sqrt{m} \left(\log
        \frac{2N}{m} \r)^{1/p}
  + 6 \max _{i\leq N} |X_i|  \r)} \\
  &  \leq &  (1+ 2 \log m ) \exp\left( -  K^p \ \sqrt{m} \log \frac{2
      N}{m} \r) .
\end{eqnarray*}
However we will not pursue this direction here.
}
\end{remark}

\section{Kannan-Lov\'asz-Simonovits question}
\label{almost_isometric_section}

In this section, we answer the question presented in the introduction:
{\em Let $K$ be an isotropic convex body in $\R^n$.  Given $\eps>0$,
  how many independent points $X_i$ uniformly distributed on $K$ are
  needed for the empirical covariance matrix to approximate the
  identity up to $\eps$ with overwhelming probability?}

Let $X\in\R^n$ be a centered random vector with covariance matrix
$\Sigma$ and consider $N$ independent random vectors $(X_i)_{i\le N}$
distributed as $X$. Using empirical processes tools, we first prove a
more general statement (Proposition \ref{almost_isometric_prop}) and
then give applications to approximation of the empirical covariance
matrix and to estimates of different norms of the matrix
$A=A^{(N)}$. In a final subsection we give a more elementary proof of
the case ($p=2$) that corresponds to the original question in
\cite{KLS}.

\subsection{Approximation of covariance matrix}
\label{cov_matrix_subsection}

First note that
because of the linear invariance, (\ref{isotropic case}) implies
$$
  \Big\|\frac{1}{N}\sum_{i=1}^N
 X_i\otimes X_i - \Sigma\Big\|\le \varepsilon \|\Sigma\|.
$$
 Therefore without loss of generality we restrict ourselves to the
 case when the covariance matrix is the identity.

\begin{thm}\label{inertia_matrix_thm}
  Let $X_1,\ldots,X_N$ be i.i.d. random vectors, distributed according
  to an isotropic, log-concave probability measure on $\R^n$. For
  every $\varepsilon \in (0,1)$ and $t \ge 1$, there exists
  $C(\varepsilon,t) > 0$, such that if $C(\varepsilon,t)n \le N$,
  then with probability at least $1 - e^{-ct\sqrt{n}}$,
\begin{align}\label{inertia_matrix_ineq}
\Big\|\frac{1}{N}\sum_{i=1}^N X_i\otimes X_i - \Id\Big\| \le \varepsilon,
\end{align}
where $c >0$ is an absolute constant. Moreover, one can take
$C(\varepsilon,t) = Ct^4\varepsilon^{-2}\log^2(2t^2\varepsilon^{-2})$,
where $C> 0$ is an absolute constant.
\end{thm}

\medskip

Since for a symmetric matrix $M$, one has $\|M\| = \sup_{y\in S^{n-1}}
\langle M y, y\rangle$ and $\E \langle X_i, y\rangle^2 = |y|^2$, one
can rewrite (\ref{inertia_matrix_ineq}) as
\begin{displaymath}
  \sup_{y\in S^{n-1}}\Big|\frac{1}{N}\sum_{i=1}^N(\langle
  X_i,y\rangle^2 - \E\langle X_i,y\rangle^2)\Big|\le \varepsilon.
\end{displaymath}

This way approximating the covariance matrix becomes a special case of
a more general problem, concerning the uniform approximation of the
moments of one dimensional marginals of an isotropic log-concave
measure by their empirical counterparts. In particular, Theorem
\ref{inertia_matrix_thm} is implied by the following result.

\begin{thm}\label{empirical_moments_thm}
  Let $X_1,\ldots,X_N$ be i.i.d. random vectors, distributed according
  to an isotropic, log-concave probability measure on $\R^n$. For any
  $p \ge 2$ and for every $\varepsilon \in (0,1)$ and $t \ge 1$, there
  exists $C(\varepsilon,t,p) > 0$, such that if
  $C(\varepsilon,t,p)n^{p/2} \le N$, then with
  probability at least $1 - e^{-c_pt\sqrt{n}}$ (where $c_p>0$ depends
  only on $p$),
\begin{align}\label{empirical_moments_ineq}
 \sup_{y\in S^{n-1}}\Big|\frac{1}{N}\sum_{i=1}^N(|\langle X_i,y\rangle|^p -
 \E|\langle X_i,y\rangle|^p)\Big|\le \varepsilon.
\end{align}
Moreover, one can take
$C(\varepsilon,t,p) = C_p t^{2p}\varepsilon^{-2}\log^{2p-2}(2
t^{2}\varepsilon^{-2})$,
where $C_p$ depends only on $p$.
\end{thm}

\medskip

\begin{remark}
  {\rm
 Proofs of both Theorems, \ref{inertia_matrix_thm} and
\ref{empirical_moments_thm}, use Theorem \ref{norm_estimate_thm} which
requires the condition $N\leq \exp(\sqrt{n})$.  For larger $N$,
however,
the result follows  by a formal argument. Assume
that  the statement has been proved  for $N\leq \exp(\sqrt{n})$ and
assume that  $N > \exp(\sqrt{n})$. Let $X_i = \{X_i(k)\} _{k=1}^n
\in \R^n$, $i\leq N$, be the random vectors under
consideration. Pick  the smallest $m$ such that $N\leq
\exp(\sqrt{m})$. Clearly, $m > n$. Now consider random vectors $Y_i
= \{Y_i(k)\} _{k=1}^m \in \R^m$, $i\leq N$, defined by $Y_i(k) =
X_i(k)$ for $k\leq n$ and $Y_i (k) = g_{ik}$ for $k>n$, where $g_{ik}$
are independent Gaussian ${\cal{N}}(0, 1)$ random variables. Then
$Y_i$'s are isotropic log-concave random vectors to which the
result can be applied.  Identifying $y=\{y(k)\} _{k=1}^n \in S^{n-1}$ with
$z=\{z(k)\} _{k=1}^m \in S^{m-1}$, defined by $z(k)=y(k)$ for $k\leq
n$, $z(k)=0$ for $k>n$, we get
\begin{eqnarray*}
  \lefteqn{
 \sup_{y\in S^{n-1}}\Big|\frac{1}{N}\sum_{i=1}^N(|\langle
 X_i,y\rangle|^p -
 \E|\langle X_i,y\rangle|^p)\Big|}\\
 &\le & \sup_{y\in S^{m-1}}\Big|\frac{1}{N}\sum_{i=1}^N(|\langle
 Y_i,y\rangle|^p
  - \E|\langle Y_i,y\rangle|^p)\Big|\le \varepsilon
\end{eqnarray*}
with probability even higher than claimed.  Thus in  the  proofs of both
theorems we may assume without loss of generality that $N\leq
\exp{(\sqrt{n})}$.
}
\end{remark}

\bigskip

In the first step of the proof of Theorem \ref{empirical_moments_thm}
we shall use some tools from the probability in Banach spaces, in
particular classical symmetrization and contraction methods as in
\cite{GR} and \cite{Men}.  These tools work for general empirical
processes and are not necessary in our setting since we are dealing
more specifically with powers of linear forms. We choose this
approach, though, as it requires less computations and leads to a
unified, simpler and more transparent presentation.

\medskip
Theorem \ref{empirical_moments_thm} is an easy consequence of the
following technical proposition applied with $s = t$.

\begin{prop}\label{almost_isometric_prop} In the setting of Theorem
  \ref{empirical_moments_thm}, if $ n \le N \le e^{\sqrt n}$, then
  for any $s,t \ge 1$,   the estimate
\begin{align}\label{almost_isometric_ineq}
\sup_{y\in S^{n-1}}&\Big|\frac{1}{N}\sum_{i=1}^N(|\langle
X_i,y\rangle|^p - \E|\langle X_i,y\rangle|^p)\Big| \nonumber\\
&\le
C^{p-1}ts^{p-1}p\log^{p-1}\Big(\frac{2N}{n}\Big)\sqrt{\frac{n}{N}} +
\frac{C^ps^pn^{p/2}}{N} + C^pp^p\Big(\frac{n}{2N}\Big)^s
\end{align}
holds with probability at least
\begin{displaymath}
  1 - \exp(-c s\sqrt{n}) -
  \exp\big(-c_p\min \{ u, v \}\big)
\end{displaymath}
where $u = t^2s^{2p-2}n\log^{2p-2}(2N/n)$,
$v= ts^{-1}\sqrt{Nn} / \log(2N/n)$, $C, c >0$ are absolute constants
and $c_p >0 $ depends on $p$ only.
\end{prop}

\begin{remark}
  {\rm
The two parameters $s$ and $t$ play different role in the proof
and reflect different asymptotic behavior of the probability with
which (\ref{almost_isometric_prop}) holds. The first parameter $s$ is
related to a level of truncation of linear forms whereas the second is
a factor in the deviation when one deals only with the truncated part.
For instance, by taking $ s = t^{1/2}$, it allows us to get a
probability converging to one as $t\to \infty$, if both dimensions are
fixed.
}
\end{remark}

\paragraph{}
Before we proceed to the proof of the above proposition, let us
introduce some tools from the classical theory of probability in
Banach spaces. Below, $\varepsilon_1,\ldots,\varepsilon_N$ will always
denote a sequence of independent Rademacher variables, independent of
the sequence $X_1,\ldots,X_N$.

\begin{lemma}[Contraction principle, see \cite{LT}, Theorem
  4.12]\label{contraction_principle}
  Let $F\colon \R^+\to \R_+$ be convex and increasing. Let further
  $\varphi_i\colon \R \to \R$, $i\le N$ be 1-Lipschitz with
  $\varphi_i(0)=0$. Then, for any bounded set $T \subset \R^N$,
\begin{displaymath}
\E F\Big(\frac{1}{2}\sup_{t\in T}\Big|\sum_{i=1}^N \varepsilon_i
\varphi_i(t_i)\Big|\Big) \le \E F\Big(\sup_{t\in T}\Big|\sum_{i=1}^N
\varepsilon_i t_i\Big|\Big).
\end{displaymath}
\end{lemma}

Using  standard symmetrization inequalities
for sums of independent random variables (see e.g., Chapter 2.3. of
\cite{vVW})
and applying the lemma with $F \equiv 1$, and $\varphi_i (s)=
\frac{|s|^p \wedge B^p}{pB^{p-1}}$ for $s \in \R$,
 we obtain the following corollary.

\begin{cor}
\label{bounded_functions}
  Let $\mathcal{F}$ be a family of functions, uniformly bounded by
  $B >0$. Then for any independent random variables $X_1,\ldots,X_N$ and
  any $p \ge 1$, we have
\begin{displaymath}
\E\sup_{f \in \mathcal{F}} \Big|\sum_{i=1}^N(|f(X_i)|^p -
\E|f(X_i)|^p)\Big| \le 4pB^{p-1}\E\sup_{f\in
  \mathcal{F}}\Big|\sum_{i=1}^N \varepsilon_i f(X_i)\Big|
\end{displaymath}
\end{cor}

We will also use the celebrated Talagrand's concentration inequality
for suprema of bounded empirical processes \cite{T}. The version
from \cite{KR}  presented below, provides the best known
constants in this inequality (we will however not take advantage of
explicit constants). For a simple proof (with worse constants) we
refer the reader to \cite{L1,L}

\begin{lemma}[\cite{KR}, Theorem 1.1]\label{Klein_Rio}
Let $X_1,X_2,\ldots,X_N$ be independent random variables with
values in a measurable space $(\mathcal{S},\mathcal{B})$ and let
$\mathcal{F}$ be a countable class of measurable functions $f
\colon \mathcal{S}\to [-a,a]$, such that for all $i$, $\E f(X_i)
=0$. Consider the random variable
\begin{displaymath}
Z =\sup_{f\in \mathcal{F}} \sum_{i=1}^N f(X_i).
\end{displaymath}
Then, for all $t \ge 0$,
\begin{displaymath}
\PP(Z \ge \E Z +t) \le \exp\Big(-\frac{t^2}{2(\sigma^2 + 2a\E Z) +
3at}\Big),
\end{displaymath}
where
\begin{displaymath}
\sigma^2 = \sup_{f\in \mathcal{F}}\sum_{i=1}^N\E f(X_i)^2.
\end{displaymath}
\end{lemma}

\noindent {\bf Proof of Proposition \ref{almost_isometric_prop}{\ \ }}
For simplicity, throughout this proof we will use the letter $C$ to
denote absolute constants, whose values may change from line to line.

For $B > 1$ (to be specified later) consider
\begin{align*}
\E&\sup_{y\in S^{n-1}}\Big| \sum_{i=1}^N \Big( (|\langle
X_i,y\rangle|\wedge B)^p - \E (|\langle X_i,y\rangle|\wedge
B)^p\Big)\Big|\\
&\le 4pB^{p-1}\E\sup_{y\in S^{n-1}}\Big| \sum_{i=1}^N \varepsilon_i
(|\langle X_i,y\rangle|\wedge B)\Big|,
\end{align*}
where the last line follows from  Corollary \ref{bounded_functions}. The
function $t\mapsto |t|\wedge B$ is a contraction, so
\begin{align*}
\E&\sup_{y\in S^{n-1}}\Big| \sum_{i=1}^N \Big( (|\langle
X_i,y\rangle|\wedge B)^p - \E (|\langle X_i,y\rangle|\wedge
B)^p\Big)\Big| \\
 &\le 8pB^{p-1}\E\sup_{y\in S^{n-1}}\Big| \sum_{i=1}^N \varepsilon_i
 \langle X_i,y\rangle\Big| \le 8pB^{p-1}\E\Big|\sum_{i=1}^N
 \varepsilon_iX_i\Big|\\
&\le 8pB^{p-1}\sqrt{Nn}.
\end{align*}
Since by (\ref{psi1}), $\E(|\langle X_i,y\rangle|\wedge B)^{2p}\le
C^{2p} p^{2p}$,
Lemma \ref{Klein_Rio} implies that for $t \ge 1$, with probability at
least
\begin{align}\label{probability_estimate}
1- \exp&\Big(-\frac{64B^{2p-2}t^2Nn}{2NC^{2p}p^{2p} + 32pB^{2p-1}
  \sqrt{Nn}+ 24pB^{2p-1}t\sqrt{Nn}}\Big) \nonumber\\
&\ge 1-  \exp(-c_p\min(t^2nB^{2p-2}, t\sqrt{Nn}/B)),
\end{align}
one has
\begin{align}\label{bounded_part}
&\sup_{y\in S^{n-1}}\Big| \sum_{i=1}^N \Big( (|\langle
X_i,y\rangle|\wedge B)^p - \E (|\langle X_i,y\rangle|\wedge
B)^p\Big)\Big| \le 16t pB^{p-1}\sqrt{Nn}.
\end{align}

Observe that
\begin{align*}
  \sup_{y\in S^{n-1}}&\Big|\frac{1}{N} \sum_{i=1}^N(|\langle
  X_i,y\rangle|^p - \E |\langle X_i,y\rangle|^p) \Big| \\
  \le& \sup_{y\in S^{n-1}} \Big| \sum_{i=1}^N \frac{1}{N}(|\langle
  X_i,y\rangle|\wedge B)^p  - \E (|\langle X_i,y\rangle|\wedge
  B)^p)\Big|\\
  &+\sup_{y\in S^{n-1}}\frac{1}{N}\sum_{i=1}^N (|\langle
  X_i,y\rangle|^p - B^p)\ind{|\langle X_i,y\rangle|\ge B} \\
  &+ \sup_{y\in S^{n-1}}\frac{1}{N}\E\sum_{i=1}^N (|\langle
  X_i,y\rangle|^p - B^p)\ind{|\langle X_i,y\rangle|\ge B},
\end{align*}

Each of the obtained three terms is estimated  separately, with the
first term already discussed in  (\ref{bounded_part})
and  (\ref{probability_estimate}).
By (\ref{psi1}) and Chebyshev's inequality we have
\begin{displaymath}
\E  |\langle X_i,y\rangle|^p \ind{|\langle X_i,y\rangle| \ge B} \le
\|\langle X_i,y\rangle\|_{2p}^{p}\sqrt{\PP(|\langle X_i,y\rangle|\ge
  B)} \le C^pp^{p}e^{-B/C}.
\end{displaymath}
 Together with the previous inequalities  this implies that
\begin{align}\label{three_term_estimate}
\sup_{y\in S^{n-1}}&\Big|\frac{1}{N} \sum_{i=1}^N(|\langle
X_i,y\rangle|^p - \E |\langle X_i,y\rangle|^p) \Big|\nonumber\\
&\le 16tpB^{p-1} \sqrt{\frac{n}{N}} +  \sup_{y\in
  S^{n-1}}\frac{1}{N}\sum_{i=1}^N |\langle
X_i,y\rangle|^p\ind{|\langle X_i,y\rangle| \ge B} + C^pp^{p}e^{-B/C},
\end{align}
with probability at least
$$
1- \exp(-c_p\min(t^2nB^{2p-2},
t\sqrt{Nn}/B)).
$$

 Thus it remains to estimate $\sup_{y\in
  S^{n-1}}\sum_{i=1}^N |\langle X_i,y\rangle|^p\ind{|\langle
  X_i,y\rangle| \ge B}$.
To this end we use Theorem \ref{norm_estimate_thm} and
Remark~\ref{rem_m-to-n}.  It follows that for $s \ge 1$, with
probability at least $1 - e^{-cs\sqrt{n}}$, we have,  for all $m \le N$ and all
$z \in S^{N-1}$ with $|\supp z| = m$,
\begin{equation}
\label{consequence}
\Big|\sum_{i=1}^N z_i X_i\Big| \le Cs\Big(\sqrt{n} +
\sqrt{m}\log\Big(\frac{2N}{n}\Big)\Big).
\end{equation}
Dualizing this estimate and using the fact that for $p \ge 2$, the
$\ell_p$ norm is dominated by the $\ell_2$ norm, we obtain,
for any set $E \subset
\{1,\ldots,N\}$,
\begin{eqnarray}
  \label{short_support_estimate-1}
\sup_{y \in S^{n-1}}\left(\sum_{i \in E} |\langle X_i,y \rangle|^p
\right)^{1/p} &\le &
\sup_{y \in S^{n-1}}\left(\sum_{i \in E} |\langle X_i,y \rangle|^2
\right)^{1/2} \nonumber\\
  &\le & Cs\Big(\sqrt{n} +
\sqrt{|E|}\log\Big(\frac{2N}{n}\Big)\Big).
\end{eqnarray}
For   an arbitrary ${y \in S^{n-1}}$
let $E_B = E_B(y) := \{i\le N \colon |\langle
X_i,y\rangle|\ge B\}$.  Then, by (\ref{short_support_estimate-1}),
$$
B |E_B|^{1/2} \le
\left(\sum_{i \in {E_B}} |\langle X_i,y \rangle|^2
\right)^{1/2}
\le Cs\Big(\sqrt{n} +
\sqrt{|E_B|}\log\Big(\frac{2N}{n}\Big)\Big).
$$
Thus, whenever
\begin{equation}
  \label{B-condition}
 B \ge 2Cs \log\Big(\frac{2N}{n}\Big),
\end{equation}
we obtain  (for a different absolute constant  $C$),
$$
|E_B| \le C s^2 n B^{-2}.
$$

This combined with (\ref{short_support_estimate-1}) implies,
after taking the $p$'th powers and  again adjusting constants,
that with probability at
least $1 - e^{-cs\sqrt{n}}$, for all $y \in S^{n-1}$,
\begin{align*}
\sum_{i =1}^N |\langle X_i,y \rangle|^p \ind{|\langle X_i,y\rangle|
  \ge B} &= \sum_{i \in E_B} |\langle X_i,y \rangle|^p\\
&\le C^ps^p\Big(n^{p/2} +
n^{p/2}s^pB^{-p}\log^p\Big(\frac{2N}{n}\Big)\Big).
\end{align*}

Setting $B = 2Cs\log(2N/n)$, so  that
(\ref{B-condition}) is satisfied,
 and combining the resulting estimate with
(\ref{three_term_estimate}), we get
\begin{align*}
\sup_{y\in S^{n-1}}&\Big|\frac{1}{N} \sum_{i=1}^N(|\langle
X_i,y\rangle|^p - \E |\langle X_i,y\rangle|^p) \Big|\nonumber\\
&\le
16C^{p-1}ts^{p-1}p\log^{p-1}\Big(\frac{2N}{n}\Big)\sqrt{\frac{n}{N}} +
\frac{C^ps^pn^{p/2}}{N} + C^pp^p\Big(\frac{n}{2N}\Big)^s,
\end{align*}
with probability at least
\begin{displaymath}
1 - \exp(-c s\sqrt{n}) -
\exp\Big(-c_p\min\Big(t^2s^{2p-2}n\log^{2p-2}(2N/n),
\frac{ts^{-1}\sqrt{Nn}}{\log(2N/n)}\Big)\Big).
\end{displaymath}
This  completes the proof of
Proposition~\ref{almost_isometric_prop},
\qed

\begin{remark}
  {\rm
 Let $G\in\R^n$ be a standard Gaussian
vector with the identity as the covariance matrix and let $h$ be a
standard Gaussian random variable. Assume that $h$ and $G$ are
independent and put $X=h\,G\in\R^n$. Clearly its covariance matrix
is the identity  and it is easy to check that $\|\langle
X,y\rangle\|_{\psi_1}\le c|y|$, for every $y\in\R^n$, where $c$ is a
numerical constant. Nevertheless, it is known from \cite{A2} that
$X$ does not satisfy the conclusion of Lemma \ref{max}; in fact the
density of $X$ is not log-concave. Now let us consider the matrix
$A=A^{(N)}$ with i.i.d. copies $X_i=h_i \,G_i$, $i=1,\dots, N$ as
columns with $N\le e^n$, where $(h_i)$ are i.i.d copies of $h$ and
similarly $(G_i)$ i.i.d copies of $G$, $(h_i)$ and $(G_i)$
independent. One can check that
\begin{eqnarray*}
\E \sup_{y\in S^{n-1}}{\frac{1} {N}}\sum_1^N  &&\!\!\!\!\!\!\!\!\!\!\!\!
 |\langle X_i,y\rangle|^2=
\E \sup_{y\in S^{n-1}}{\frac{1}{N}}\sum_1^N h_i^2|\langle G_i,y\rangle
|^2\\
& \ge &
\E \sup_{i}{\frac{1}{ N}} h_i^2|G_i|^2\ge c{\frac{n} {N} }\log N
\end{eqnarray*}
where $c>0$ is a numerical constant. Thus $\|A\|\ge \sqrt{c{n }\log
N}$. This example shows that the sub-exponential decay of linear
forms ($\psi_1$ norm bounded) is not sufficient for our problem.
}
\end{remark}

\medskip

\begin{remark}
  {\rm
 In comparison, a sub-gaussian decay of
linear forms is sufficient. Indeed, it is known (see for instance
\cite{MPT}) that if there exists $c>0$ such that $\E \exp \left(|c
\langle X,y\rangle|^2\right)\le 2$ for every $y\in S^{n-1}$, then
(\ref{isotropic case}) holds with probability larger than
$1-\exp(-c'n)$ for some numerical constant $c'>0$.
}
\end{remark}

\medskip

\begin{remark}
  {\rm
Another non necessarily log-concave
example for which the conclusion of Theorems \ref{norm_estimate_thm}
and \ref{inertia_matrix_thm} are valid is obtained  when $\|\langle
X,y\rangle\|_{\psi_1}\le c|y|$, for every $y\in\R^n$ and $|X|\le
C\sqrt n$ where $c,C>0$ are numerical constants.
}
\end{remark}

\subsection{Additional observations}
\label{additional_obs_subsection}

We note several observations for norms of random matrices from
$\ell_2$ to $\ell_p$, $p \ne 2$.

\begin{cor}\label{ell_p_cor}
  For $1\le N \le e^{\sqrt{n}}$ let $\Gamma$ be a random $N\times n$
  matrix with rows $X_1,\ldots,X_N$. Then for $p \ge 2$, with
  probability at least $1 - e^{-c_p\sqrt{n}}$ (where $c_p >0$ depends
  only on $p$),
\begin{align}\label{ell_p_estimate}
\|\Gamma\|_{\ell_2\to \ell_p} \le C_p(N^{1/p} + n^{1/2}),
\end{align}
with $C_p > 0$ depending only on $p$.
Moreover
\begin{align}\label{ell_p_exp_estimate}
  \tilde{c}_pN^{1/p} + c\sqrt{n} \le \E \|\Gamma\|_{\ell_2\to \ell_p}
  \le \tilde{C}_p(N^{1/p} + n^{1/2}),
\end{align}
where $\tilde{C}_p,\tilde{c}_p > 0$ depend only on $p$ and $c>0$ is an
absolute constant.
\end{cor}

\proof Inequality (\ref{ell_p_estimate}) for $N \le n$ follows from
Theorem \ref{norm_estimate_thm} and the comparison between $\ell_p$
norms. For $N \ge n$, the inequality follows from Proposition
\ref{almost_isometric_prop}.

Since by log-concavity, moments and quantiles of $\|\Gamma\|_{\ell_2\to
  \ell_p}$ are equivalent, (\ref{ell_p_estimate}) implies that
\begin{displaymath}
\E \|\Gamma\|_{\ell_2\to \ell_p} \le \tilde{C}_p(N^{1/p}  + n^{1/2}).
\end{displaymath}

On the other hand, a single row of $\Gamma$ has expected Euclidean
norm of the order of $\sqrt{n}$ and a single column of $\Gamma$ has
expected $\|\cdot\|_p$ norm of the order of $c(p)N^{1/p}$, so the left
hand side of (\ref{ell_p_exp_estimate}) follows trivially.  \qed

\begin{cor}\label{ell_p_cor_1}
  For $1\le N \le e^{\sqrt{n}}$ let $\Gamma$ be a random $N\times n$
  matrix with rows $X_1,\ldots,X_N$. Then for $p \in [1, 2)$, with
  probability at least $1 - e^{-c\sqrt{n}}$ (where $c>0$ is an
  absolute constant),
\begin{align}\label{ell_p_estimate_1}
\|\Gamma\|_{\ell_2\to \ell_p} \le C(N^{1/p} + N^{1/p-1/2}n^{1/2})
\end{align}
for some absolute constant $C> 0$.
Moreover
\begin{align}\label{ell_p_exp_estimate_1}
  \tilde{c}(N^{1/p} + N^{1/p-1/2}n^{1/2}) \le \E \|\Gamma\|_{\ell_2\to
    \ell_p} \le \tilde{C}(N^{1/p} + N^{1/p - 1/2}n^{1/2}),
\end{align}
where $\tilde{C},\tilde{c} > 0$ are absolute constants.
\end{cor}

\proof
Inequality (\ref{ell_p_estimate_1}) and the right-hand side of
(\ref{ell_p_exp_estimate_1}) follow from the corresponding results for
$p = 2$, since
\begin{displaymath}
\|\Gamma\|_{\ell_2 \to \ell_p} \le N^{1/p - 1/2}\|\Gamma\|_{\ell_2\to \ell_2}.
\end{displaymath}
To prove the left-hand side of (\ref{ell_p_exp_estimate_1}), it is
enough to notice that if $1/p^\ast + 1/p = 1$, then
\begin{displaymath}
  \E \|\Gamma\|_{\ell_2\to \ell_p} \ge \E \Big|\sum_{i=1}^N
  \frac{1}{N^{1/p^\ast}}X_i\Big| \ge \tilde{c} N^{1/2-1/p^\ast}n^{1/2}
  = \tilde{c}N^{1/p-1/2}n^{1/2}
\end{displaymath}
and the expected $\ell_p$ norm of a single column of $\Gamma$ is at
least $\tilde{c}N^{1/p}$.
\qed

One can also obtain an almost-isometric result for $p \in [1,2)$.

\begin{thm}\label{empirical_moments_thm12}
  Let $X_1,\ldots,X_N$ be i.i.d. random vectors, distributed according
  to an isotropic, log-concave probability measure on $\R^n$. For any
  $p \in [1,2)$ and for every $\varepsilon \in (0,1)$ and $t \ge 1$,
  there exists $C(\varepsilon,t)>0$, such that if $C(\varepsilon)n \le
  N\le e^{\sqrt{n}}$, then with probability at least $1 -
  e^{-ct\sqrt{n}}$ (where $c > 0$ is an absolute constant),
\begin{align}\label{empirical_moments_ineq12}
\sup_{y\in S^{n-1}}\Big|\frac{1}{N}\sum_{i=1}^N(|\langle X_i,y\rangle|^p
- \E|\langle X_i,y\rangle|^p)\Big|\le \varepsilon.
\end{align}
Moreover, one can take $C(\varepsilon,t) =
Ct^{2p}\varepsilon^{-2}\log^{2p-2}(2t^{2p}\varepsilon^{-2})$, where
$C>0$ is an absolute constant.

\end{thm}

\proof Since the proof differs only by technical details from the
corresponding argument for $p \ge 2$, we will just indicate the
necessary changes. We will use the notation from the proof of
Proposition \ref{almost_isometric_prop}.

Just as before, we truncate at the level of $Ct\log(2N/n)$ and use the
contraction principle to handle the bounded part of the process. As
for the unbounded part, we also proceed as before, however now we use
the comparison between the $\ell_2^k$ and $\ell_p^k$ norm for $p < 2$
and $k = |E_B| \le n$, which yields

\begin{align*}
\sup_{y\in S^{n-1}}&\Big|\frac{1}{N} \sum_{i=1}^N(|\langle
X_i,y\rangle|^p - \E |\langle X_i,y\rangle|^p) \Big|\nonumber\\
&\le 16C^{p-1}t^pp\log\Big(\frac{2N}{n}\Big)^{p-1}\sqrt{\frac{n}{N}} +
\frac{C^pt^pn}{N} + \frac{C^pp^pn}{N},
\end{align*}
with probability at least
\begin{displaymath}
1 - \exp(-ct\sqrt{n}) -
\exp(-c\min(t^2n\log^{2p-2}(2N/n),\sqrt{Nn}/\log(2N/n)))
\end{displaymath}
(the constants in the exponents can be made independent of $p$, since
now $p$ runs over a bounded interval). This allows us to finish the
proof.
\qed

\begin{remark}
  {\rm
 The isomorphic result for $p=1$ was proven
in \cite{GM}. The  same paper also  considers  $p \in (0,1)$.
}
\end{remark}

\subsection{Elementary approach for $p=2$}
\label{elementary_section}

As announced earlier  we will now briefly describe
a more elementary proof of Theorem \ref{inertia_matrix_thm} and
Theorem \ref{empirical_moments_thm} for $p=2$.
In this case, the classical Bernstein inequality and a net argument on
the sphere may replace the contraction principle and concentration of
measure for empirical processes, that have been used -- via Lemma
\ref{Klein_Rio} -- to prove (\ref{bounded_part}).  The remaining part
of the proof is left unchanged.

The key point is the following well known observation:
\begin{lemma} \label{from_net_to_sphere}
Let $x_i$, $i = 1,2,\ldots,N$, be arbitrary vectors in $\R^n$.
Let $\eps \in (0,1)$ and let
$ \mathcal{N} $ be a $c\eps$-net of
$S^{n-1}$, for some constant $c\in (0,1)$. If we have
$$
\sup_{y\in \mathcal{N}}\Big|\frac{1}{N}\sum_{i=1}^N(\langle x_i,y\rangle^2 -
1)\Big|
\le \eps$$
 then
$$
\sup_{y\in S^{n-1}}\Big|\frac{1}{N}\sum_{i=1}^N(\langle x_i,y\rangle^2 -
1)\Big|
\le c'\eps$$
where $c'$ depends on $c$.
\end{lemma}

We postpone the proof of this Lemma and pass to the proof of
 Theorems \ref{inertia_matrix_thm} and
\ref{empirical_moments_thm}.

Fix a $c\varepsilon$-net $\mathcal{N}$ of $S^{n-1}$ of cardinality at
most $(3/c\varepsilon)^n$, and $B >0$ to be determined later.  Pick an
arbitrary $y \in S^{n-1}$.

For the reader's convenience recall Bernstein's inequality.

\begin{prop}
[Bernstein's inequality, cf.~e.g.,~\cite{vVW}]
\label{Bernstein}
Let $Z_i$ be independent random variables, centered and such that
$|Z_i|\le a$ for all $1\le i\le N$.
Put
$Z ={\frac{1}{ N}} \sum_{i=1}^N Z_i$.
Then for all $\tau \ge 0$,
\begin{displaymath}
\PP(Z \ge \tau) \le \exp\Big(-\frac{\tau^2N}{2(\sigma^2 + a\tau/3)}\Big),
\end{displaymath}
where
\begin{displaymath}
\sigma^2 = (1/N)\sum_{i=1}^N Var(Z_i).
\end{displaymath}
\end{prop}

In our case $Z_i =
(|\langle X_i,y\rangle|\wedge B)^2 - \E (|\langle X_i,y\rangle|\wedge
B)^2 $, for $1 \le i \le N$, $a = B^2$.
Since $\E (|\langle X_i,y\rangle|)^2 =1$ then
 (\ref{psi1}) implies
$$
Var(Z_i) \le \E (|\langle X_i,y\rangle|\wedge B)^4 \le c.
$$
  Setting  $\tau = t B
\sqrt{n/N}$ we infer  that
$$
\Big| \frac{1}{N}\sum_{i=1}^N \Big( (|\langle X_i,y\rangle|\wedge B)^2
- \E (|\langle X_i,y\rangle|\wedge B)^2 \Big)\Big| \ge t B \sqrt{n/N}
$$
with probability at most
$$
 \exp \Bigl(- c \min \bigl(t^2 B^2 n, t \sqrt{Nn}/B \bigr)\Bigr).
$$
By  the union bound,
\begin{align}\label{bounded_part_eps_net}
&\sup_{y\in \mathcal{N}}
\Big|\frac{1}{N} \sum_{i=1}^N \Big( (|\langle X_i,y\rangle|\wedge B)^2 - \E
(|\langle X_i,y\rangle|\wedge B)^2\Big)\Big| \le t B \sqrt{n/N},
\end{align}
with probability at least
\begin{displaymath}
1-  \exp\Big(n\log\Big(\frac{3}{c\varepsilon}\Big) -c\min(t^2nB^2,
 t\sqrt{Nn}/B)\Big).
\end{displaymath}
This estimate corresponds to (\ref{bounded_part}).

Using this estimate with $B = Ct\log(2N/n)$ and handling the unbounded
part  the same way as in Proposition \ref{almost_isometric_prop}
(see the  argument that follows  (\ref{bounded_part}))
we obtain
\begin{align}\label{explicit_bound_on_net}
  \sup_{y\in \mathcal{N}}&\Big|\frac{1}{N} \sum_{i=1}^N(|\langle
  X_i,y\rangle|^2
 - \E |\langle X_i,y\rangle|^2) \Big|\nonumber\\
  &\le C t^{2}\log\Big(\frac{2N}{n}\Big) \sqrt{\frac{n}{N}}
  + \frac{C^2 t^2 n}{N} + \frac{4 C^2 n}{N},
\end{align}
with probability at least
\begin{displaymath}
  1 - \exp(-ct\sqrt{n}) -
  \exp\Big(n\log\Big(\frac{3}{c\varepsilon}\Big) -
  c\min\Big(t^{4}n\log^{2}(2N/n),
\frac{\sqrt{Nn}}{C\log(2N/n)}\Big)\Big).
\end{displaymath}
This corresponds to the  estimates in Proposition
\ref{almost_isometric_prop}  (for  $s=t$).

\bigskip

Now, for $N \ge C (\ep, t) n $, and $C(\ep, t)$ sufficiently large,
the right hand side of (\ref{explicit_bound_on_net}) is at most
$\varepsilon$ and  $5/c\varepsilon \le 2N/n$ which leads to the
probability above to be at least  $1 - \exp(-c t\sqrt{n})$.
So with the same probability we get
\begin{displaymath}
  \sup_{y\in \mathcal{N}}\Big|\frac{1}{N} \sum_{i=1}^N(|\langle
  X_i,y\rangle|^2
 - \E |\langle X_i,y\rangle|^2) \Big| \le \varepsilon.
\end{displaymath}
We can now conclude by Lemma \ref{from_net_to_sphere} applied
pointwise with $x_i = X_i(\omega)$ for $\omega$ from the event on
which our estimates hold (recall that by the isotropicity assumption
we have $\E |\langle X_i,y\rangle|^2=1$).

\medskip

\noindent {\bf Proof of Lemma \ref{from_net_to_sphere}}
Consider the semi-norm $\|\cdot\|$ on $\R^n$  defined by
$$
\|y\|= \Big(\frac{1}{N} \sum_{i=1}^N|\langle
  x_i,y\rangle|^2  \Big)^{1/2},
$$
for $y \in \R^n$. Our assumptions imply that
$$
 1 - \eps \le \sqrt{1-\eps} \le
 \sup_{y\in \mathcal{N}} \|y\| \le \sqrt{1+\eps}\le 1 + \eps/2.
$$
The triangle inequality and homogeneity of $\|\cdot\|$ imply, by a
standard argument, that
$$
   \sup_{y\in S^{n-1}} \|y\| \le (1 + \eps/2)(1-c\eps)^{-1} \le 1 +
   \delta,
$$
where
\begin{displaymath}
\delta = \frac{1+5c-3c^2}{2(1-c)}\varepsilon
\end{displaymath}

To get a lower estimate, write an arbitrary $y \in S^{n-1}$ in the
form $y = y_1 + c\eps y_2$, with $y_1 \in \mathcal{N}$ and $y_2 \in
S^{n-1}$.
Then $\|y\|\ge \|y_1\| - c\eps \|y_2\| \ge (1-\eps) - c \eps
(1+\delta)
\ge 1 - \delta_1$,
where
\begin{displaymath}
\delta_1 = \frac{2+c+3c^2-3c^3}{2(1-c)}\varepsilon.
\end{displaymath}
Thus for all $y\in S^{n-1}$, $|\|y\|-1| \le c_1\varepsilon$ for some $c_1$ depending only on $c$. In particular $\|y\|\in [0,1+c_1]$. Using the fact that the function $t\mapsto t^2$ is Lipschitz with constant $2(1+c_1)$ on the interval $[0,1+c_1]$, we conclude that
\begin{displaymath}
\sup_{y\in S^{n-1}}\Big|\frac{1}{N}\sum_{i=1}^N(\langle x_i,y\rangle^2 -
1)\Big|
\le c'\eps,
\end{displaymath}
where $c' = 2c_1(1+c_1)$ depends only on $c$.
\qed

\vspace{1cm}

\address

\end{document}